\documentclass[11pt]{amsart}

\usepackage[expansion=false]{microtype}
\usepackage[a4paper,margin=28mm]{geometry}
\usepackage{amsmath,amssymb,amsthm,mathtools}
\usepackage{enumitem}
\usepackage{booktabs}
\usepackage{xcolor}
\usepackage{hyperref}
\usepackage{parskip}
\hypersetup{
  colorlinks=true,
  linkcolor=blue!50!black,
  citecolor=blue!50!black,
  urlcolor=blue!50!black,
  pdftitle={Decorated broken ribbons and a Murnaghan--Nakayama rule for
            dual canonical Grothendieck functions}
}

\DeclareRobustCommand{\S}{\ensuremath{\mathsection}}

\numberwithin{equation}{section}

\newtheorem{theorem}{Theorem}[section]
\newtheorem{proposition}[theorem]{Proposition}
\newtheorem{lemma}[theorem]{Lemma}
\newtheorem{corollary}[theorem]{Corollary}
\newtheorem{maintheorem}{Theorem}

\theoremstyle{definition}
\newtheorem{definition}[theorem]{Definition}
\newtheorem{example}[theorem]{Example}
\newtheorem{notation}[theorem]{Notation}
\theoremstyle{remark}
\newtheorem{remark}[theorem]{Remark}

\newcommand{\Lam}{\Lambda}
\newcommand{\hLam}{\widehat\Lambda}
\newcommand{\Par}{\mathcal P}
\newcommand{\SSYT}{\operatorname{SSYT}}
\newcommand{\wt}{\operatorname{wt}}
\newcommand{\swt}{\operatorname{swt}}
\newcommand{\leg}{\operatorname{leg}}
\newcommand{\supp}{\operatorname{supp}}
\newcommand{\Inc}{\operatorname{Inc}}
\newcommand{\Dec}{\operatorname{Dec}}
\newcommand{\rows}{\operatorname{rows}}
\newcommand{\cols}{\operatorname{cols}}
\newcommand{\NW}{\operatorname{NW}}
\newcommand{\hgt}{\operatorname{ht}}
\newcommand{\Add}{\operatorname{Add}}
\newcommand{\enD}{\operatorname{end}}
\newcommand{\ev}{\operatorname{ev}}
\newcommand{\bfU}{\mathbf U}
\newcommand{\bbeta}{\boldsymbol\beta}
\newcommand{\balpha}{\boldsymbol\alpha}
\newcommand{\Hall}[2]{\left\langle #1,#2\right\rangle}
\newcommand{\ind}{\mathbf 1}
\newcommand{\sk}{/\!\!/}
\newcommand{\one}{\varnothing}
\newcommand{\sG}{\mathsf G}
\newcommand{\sg}{\mathsf g}
\newcommand{\cC}{\mathcal C}
\newcommand{\cI}{\mathcal I}
\newcommand{\cH}{\mathcal H}
\newcommand{\cB}{\mathcal B}
\newcommand{\QQ}{\mathbb Q}
\newcommand{\ZZ}{\mathbb Z}
\newcommand{\DHist}{\operatorname{DH}}

\title[Murnaghan--Nakayama rules for dual canonical Grothendieck functions]
{Hook histories and a two-parameter\\ Murnaghan--Nakayama rule for\\
dual canonical Grothendieck functions}

\author{Siddheswar Kundu, Santosha Pattanayak}

\address{
Department of Mathematics,
Indian Institute of Technology Kanpur,
Kanpur, Uttar Pradesh 208016, India,
Email: kundusidhu96@gmail.com, siddheswark@iitk.ac.in
}

\address{
Department of Mathematics,
Indian Institute of Technology Kanpur,
Kanpur, Uttar Pradesh 208016, India,
Email: santosha@iitk.ac.in
}

\subjclass[2020]{Primary 05E05; Secondary 05E10, 14M15, 19E08}
\keywords{Murnaghan--Nakayama rule, canonical stable Grothendieck
functions, dual Grothendieck functions, noncommutative Schur functions,
free fermions, $K$-homology of Grassmannians}

\date{\today}

\begin{document}

\begin{abstract}
We establish a two-parameter Murnaghan--Nakayama rule for power-sum
multiplication in the dual canonical stable Grothendieck basis.
The coefficients admit a closed boundary product formula and a
signed combinatorial description by hook histories. We determine
their exact support and prove an exact divisibility law for
off-diagonal coefficients governed by the edge-connected components
of the indexing skew diagram. We also obtain a hook expansion of
power sums, an explicit binomial reduction to the dual stable rule,
and a Grassmannian $K$-homological interpretation of the one-parameter
specialization. Finally, we deduce complementary rules for power-sum
multiplication in the canonical basis and the corresponding lowering
operation in the dual canonical basis.
\end{abstract}

\maketitle

\section{Introduction}\label{sec:introduction}

Murnaghan--Nakayama rules describe multiplication by the power sums $p_k$
in a distinguished basis of the ring of symmetric functions. The classical
rule concerns the Schur basis, and through the characteristic map it
computes the irreducible characters of the symmetric groups
\cite[Ch.~I, \S7]{Macdonald}: 
\begin{equation}\label{eq:classical-MN}
 p_k s_\lambda=\sum_{\mu}(-1)^{\hgt(\mu/\lambda)}s_\mu ,
\end{equation}
where the sum runs over all partitions $\mu\supseteq\lambda$ such that
$\mu/\lambda$ is a ribbon with $k$ boxes, and $\hgt(\mu/\lambda)$ is the
number of occupied rows minus one \cite[Ch.~I, \S3, Ex.~11]{Macdonald}.

Stable Grothendieck functions $G_\lambda$ represent Schubert classes in
the $K$-theory of Grassmannians \cite{Buch}, and their Hall duals
$g_\lambda$ represent classes in $K$-homology \cite{LamPylyavskyy}.
Yeliussizov \cite{Yeliussizov} introduced a two-parameter family
$G^{(\alpha,\beta)}_\lambda$ of \emph{canonical stable Grothendieck
functions}, which interpolates between the stable Grothendieck functions
($\alpha=0$) and their images under the involution $\omega$ ($\beta=0$), and
he studied the dual family. We normalize the duality as
\begin{equation}\label{eq:duality-intro}
 \Hall{G^{(-\alpha,-\beta)}_\rho}{g^{(\alpha,\beta)}_\lambda}=\delta_{\rho\lambda},
\end{equation}
so that $g^{(0,0)}_\lambda=s_\lambda$ and $g^{(0,1)}_\lambda$ is the dual
stable Grothendieck function of Lam and Pylyavskyy; see
Section~\ref{subsec:canonical} for the sign conventions.

There are two natural multiplication problems associated with the
Hall-dual bases
$\{G^{(-\alpha,-\beta)}_\lambda\}$ and
$\{g^{(\alpha,\beta)}_\lambda\}$.
The first concerns multiplication in the canonical basis, whereas
the second concerns multiplication in the dual canonical basis.
Nguyen, Hiep, Son and Thuy \cite{NHST} established a
Murnaghan--Nakayama rule for multiplication by ordinary power sums
in the one-parameter stable Grothendieck basis. 
The first author \cite{Kundu} obtained a two-parameter
canonical-basis rule with the deformed multipliers
$p_k(x/(1-\alpha x))$.
By Hall adjointness, these results yield corresponding lowering
rules on the dual basis.

For the dual stable Grothendieck basis, Hiep and Thuy
\cite[Thm.~2.1]{HiepThuy} recently obtained a closed
Murnaghan--Nakayama rule, which in our notation corresponds to
$g_\lambda=g^{(0,1)}_\lambda$.
Their off-diagonal coefficients are signed binomial coefficients
supported on corner-connected broken ribbons. They depend on the
number of boxes, the number of occupied rows, the number of
edge-connected components of $\mu/\lambda$, and a boundary
statistic of the pair $(\lambda,\mu)$.

In this paper, we treat the full two-parameter family
$\{g^{(\alpha,\beta)}_\lambda\}$, with $\alpha$ and $\beta$
independent. Our boundary product formula specializes to the
one-parameter description and admits a direct Schur-operator
proof, together with a signed interpretation in terms of
decorated histories. Our rule (Theorem~\ref{thm:A}) specializes at
$(\alpha,\beta)=(0,1)$ to the formula of Hiep and Thuy
(Corollary~\ref{cor:compact-binomial}).
In this specialization, their boundary statistic agrees with
our $B$, so their formula is recovered from ours.
The expansion of $p_k$ in the dual canonical basis
(Theorem~\ref{thm:C}) also specializes at $(0,1)$ to
\cite[Cor.~2.5]{HiepThuy}.

In addition, our methods yield the $K$-homological interpretation,
rational formulas for the deformed power sums, and, via Kundu's
rule, multiplication by $p_k$ in the canonical basis and the
corresponding lowering rule in the dual basis
(Theorem~\ref{thm:D}).

The multiplication problem studied here is
\begin{equation}\label{eq:raising-problem}
 p_k\,g^{(\alpha,\beta)}_\lambda
 =\sum_\mu b^{(k)}_{\lambda\mu}(\alpha,\beta)\,
              g^{(\alpha,\beta)}_\mu.
\end{equation}
We call this the \emph{raising problem}.
Since multiplication by $p_k$ is adjoint to
$p_k^\perp=k\,\partial/\partial p_k$, its Hall-dual formulation is
\[
 p_k^\perp G^{(-\alpha,-\beta)}_\mu
 =\sum_\lambda b^{(k)}_{\lambda\mu}(\alpha,\beta)\,
                 G^{(-\alpha,-\beta)}_\lambda;
\]
see Corollary~\ref{cor:G-lowering}.

Our main result determines the exact support and the coefficients
in~\eqref{eq:raising-problem}. Off the diagonal, the support consists
of corner-connected broken ribbons satisfying a bound involving
their size and number of edge-connected components. The coefficients
depend on the skew diagram together with two endpoint contacts with
the inner partition. In particular, the skew diagram alone does not
determine them.

We start from the Schur-operator theorem of Iwao, Motegi, and
Scrimshaw \cite[Thm.~3.2]{IMS25}.
We justify its specialization to constant parameters and combine
it with the hook expansion 
\begin{equation}\label{eq:hook-expansion}
 p_k=\sum_{a=0}^{k-1}(-1)^a s_{(k-a,1^a)}.
\end{equation}
This gives a finite alternating formula in terms of hook-word
histories. The cancellation of disconnected word supports is the
Fomin--Greene argument \cite[Thm.~1.3]{FominGreene}, also used in
\cite{BSZ,NguyenKk}. We then account for the weighted stays and the
supersymmetric correction, remove the auxiliary cutoff, and evaluate
the remaining contributions. The resulting boundary formula gives
the sharp support condition and the explicit coefficients.

\subsection{Main results}
We first present the closed coefficient formula for the
Murnaghan--Nakayama rule, then give its combinatorial
interpretation in terms of hook histories, and finally
present two further consequences. Throughout, we abbreviate $g_\lambda=g^{(\alpha,\beta)}_\lambda$,
append zero parts to partitions, and put $\nu_0=\infty$ for
every partition $\nu$. For a condition $C$ we write $\ind_{\{C\}}=1$ if
$C$ holds and $\ind_{\{C\}}=0$ otherwise; $d(\lambda)$ is the number of
distinct positive parts of $\lambda$; and $[z^m]F(z)$ denotes the coefficient
of $z^m$ in a power series $F(z)$. We use the convention $[z^m]F(z)=0$ for $m<0$. A \emph{broken ribbon} is a skew diagram
that contains no $2\times2$ block of boxes, without any connectivity
requirement.

Let $\lambda\subsetneq\mu$. We call $\mu/\lambda$ a \emph{corner-connected
broken ribbon} if its occupied rows form an interval $[a,b]$ and
\begin{equation}\label{eq:cc-intro}
 \mu_i\in\{\lambda_{i-1},\lambda_{i-1}+1\}\qquad(a<i\leq b);
\end{equation}
equivalently, $\mu/\lambda$ contains no $2\times2$ block and is connected
when boxes sharing a corner are also considered adjacent. For such a
diagram let $D=|\mu/\lambda|$, let $r=b-a+1$ be its number of rows, let $c$
be its number of edge-connected components, and put
\[
 \varepsilon_{\mathrm{top}}=\ind_{\{a>1,\ \mu_a=\lambda_{a-1}\}},\qquad
 \varepsilon_{\mathrm{bot}}=\ind_{\{\lambda_b=\lambda_{b+1}\}},\qquad
 B=r-1+\varepsilon_{\mathrm{top}}+\varepsilon_{\mathrm{bot}},\qquad
 A=D+c-B .
\]

\begin{maintheorem}[
Theorem~\ref{thm:compact-boundary} and
Corollaries~\ref{cor:sharp-support}
and~\ref{cor:exact-component-divisibility}]
\label{thm:A}
For every partition $\lambda$ and every $k\geq1$,
\[
 p_k\,g_\lambda=d(\lambda)\bigl((-\alpha)^k-\beta^k\bigr)\,g_\lambda
 +\sum_\mu(-1)^{r-c}(\alpha+\beta)^{c-1}\,
 [z^{\,k-D-c+1}]\,\frac{1}{(1+\alpha z)^A(1-\beta z)^B}\;g_\mu ,
\]
where $\mu$ runs over the partitions $\mu\supsetneq\lambda$ such that
$\mu/\lambda$ is a corner-connected broken ribbon with
$|\mu/\lambda|+c-1\leq k$, and $D,r,c,A,B$ are formed with respect to
$(\lambda,\mu)$. Each coefficient in the sum is a nonzero polynomial in
$\ZZ[\alpha,\beta]$, homogeneous of degree $k-|\mu/\lambda|$,
and divisible by $(\alpha+\beta)^{c-1}$ but not by
$(\alpha+\beta)^c$.
\end{maintheorem}

For $k=1$ the rule reads
$p_1g_\lambda=\sum_\mu g_\mu-(\alpha+\beta)\,d(\lambda)\,g_\lambda$, the sum
over all $\mu$ obtained from $\lambda$ by adding one box; this is the case
$k=1$ of Yeliussizov's Pieri rule \cite[Prop.~8.12]{Yeliussizov}. For
$|\mu/\lambda|=k$ only ribbons ($c=1$) occur, with coefficient
$(-1)^{\hgt(\mu/\lambda)}$, as in~\eqref{eq:classical-MN}.

Theorem~\ref{thm:A} is derived from a signed combinatorial rule, which we
now describe. The operator $U_i$ acts on a partition $\nu$ in one of the
following ways:
\begin{center}
\renewcommand{\arraystretch}{1.15}
\begin{tabular}{lll}
\toprule
condition & possible moves & weights\\
\midrule
$0<\nu_i<\nu_{i-1}$ & add a box in row $i$, or stay & $1$, $-\alpha$\\
$0=\nu_i<\nu_{i-1}$ & add a box in row $i$ & $1$\\
$\nu_i=\nu_{i-1}$ & stay & $\beta$\\
\bottomrule
\end{tabular}
\end{center}
A \emph{hook word} of length $k$ is a word
\[
 w=d_a\cdots d_1\,c_1c_2\cdots c_t,
 \qquad d_a>\cdots>d_1>c_1\leq c_2\leq\cdots\leq c_t,
 \qquad a+t=k,\ t\geq1,
\]
of positive integers; these are exactly the column reading words of
semistandard tableaux of hook shape $(t,1^a)$. We put $\leg(w)=a$. The
letters act from right to left, the letter $i$ acting by $U_i$. A
\emph{history} $h$ of $w$ from $\lambda$ is a choice of an admissible move
for each letter; 
we write $\#_\alpha(h)$ and $\#_\beta(h)$ for its numbers of
stays of weight $-\alpha$ and $\beta$. Let $\cC^+_k(\lambda,\mu)$ be the set
of pairs $(w,h)$ in which $w$ is a hook word of length $k$ whose set of
letters is an interval of integers, and $h$ is a history of $w$ from
$\lambda$ to $\mu$ containing at least one addition.

\begin{maintheorem}[Theorem~\ref{thm:connected-rule}
and Proposition~\ref{prop:broken-ribbon}]\label{thm:B}
For all partitions $\lambda\subseteq\mu$ and all $k\geq1$, the coefficient
of $g^{(\alpha,\beta)}_\mu$ in $p_kg^{(\alpha,\beta)}_\lambda$ is
\[
 \sum_{(w,h)\in\cC^+_k(\lambda,\mu)}
  (-1)^{\leg(w)}(-\alpha)^{\#_\alpha(h)}\beta^{\#_\beta(h)}
 +\delta_{\lambda\mu}\,d(\lambda)\bigl((-\alpha)^k-\beta^k\bigr).
\]
The set $\cC^+_k(\lambda,\mu)$ is finite, involves only rows
$\leq\ell(\lambda)+k$, and is empty unless $\mu/\lambda$ is a broken ribbon.
\end{maintheorem}

The boxes added during a history in $\cC^+_k(\lambda,\mu)$
form the broken ribbon $\mu/\lambda$, and its stays may be recorded as \emph{ghosts} of weight
$-\alpha$ or $\beta$ on the boundary of the intermediate partitions
(Section~\ref{sec:cancellation}). Summing over the ghosts row by row gives an interval formula
(Theorem~\ref{thm:boundary}), which collapses to the single product of
Theorem~\ref{thm:A} (Section~\ref{sec:boundary}). 

Taking $\lambda=\one$ gives a closed expansion of $p_k$ in the dual
canonical basis.

\begin{maintheorem}[Theorem~\ref{thm:empty-partition}]\label{thm:C}
For every $k\geq1$,
\[
 p_k=\sum_{\substack{m\geq1,\ j\geq0\\ m+j\leq k}}
 (-1)^j\,[t^{\,k-m-j}]\,\frac{1}{(1+\alpha t)^m(1-\beta t)^{j+1}}\;
 g^{(\alpha,\beta)}_{(m,1^j)} .
\]
\end{maintheorem}

Finally, Kundu's rule for deformed power sums yields, after inverting the
deformation, a rule for multiplication by $p_k$ in the canonical basis and,
by Hall duality, a lowering rule for the dual basis.

\begin{maintheorem}[Proposition~\ref{prop:G-side} and
Theorem~\ref{thm:lowering}]\label{thm:D}
There are explicit coefficients $L^{(k)}_{\lambda\mu}(\alpha,\beta)$
(see~\eqref{eq:L-def}), which vanish unless $\mu/\lambda$ is connected with
$|\mu/\lambda|\geq k$, and otherwise depend only on the size and the numbers
of rows and columns of $\mu/\lambda$, such that
\[
 p_k\,G^{(-\alpha,-\beta)}_\lambda=\sum_\nu L^{(k)}_{\lambda\nu}(\alpha,\beta)\,
 G^{(-\alpha,-\beta)}_\nu,
 \qquad
 p_k^\perp g^{(\alpha,\beta)}_\mu=\sum_\lambda L^{(k)}_{\lambda\mu}(\alpha,\beta)\,
 g^{(\alpha,\beta)}_\lambda .
\]
\end{maintheorem}

The raising coefficients behave very differently from the lowering ones. The terms with $|\mu/\lambda|=k$ recover the classical ribbon
coefficients (Corollary~\ref{cor:top}). Terms with fewer than
$k$ added boxes can have several edge-connected components. For instance,
\begin{equation}\label{eq:intro-disconnected}
 [g_{(3,2)}]\,p_3g_{(2,1)}=\alpha+\beta ,
\end{equation}
although $(3,2)/(2,1)$ consists of two boxes that touch only at a corner
(Example~\ref{ex:disconnected}). Moreover, the coefficients are not
functions of the skew diagram. The diagonal coefficient
$d(\lambda)\bigl((-\alpha)^k-\beta^k\bigr)$ already depends on $\lambda$, and
off the diagonal
\begin{equation}\label{eq:intro-counterexample}
 b^{(2)}_{(2,1),(2,2)}=\beta-\alpha,
 \qquad
 b^{(2)}_{(3,1),(3,2)}=-2\alpha ,
\end{equation}
although both skew diagrams are the single box $(2,2)$ (Example~\ref{prop:not-shape-only}). This contrasts with the $G$-side
rules of \cite{NHST,Kundu} and with the lowering rule of
Theorem~\ref{thm:D}. Theorem~\ref{thm:A} shows that, beyond the skew diagram $\mu/\lambda$, the
coefficients depend only on the two endpoint contacts
$\varepsilon_{\mathrm{top}}$ and $\varepsilon_{\mathrm{bot}}$; the examples
above show that the skew diagram alone does not suffice.

\subsection{Organization}

Section~\ref{sec:preliminaries} fixes conventions and proves the finite
specialization of the Schur-operator theorem.
Section~\ref{sec:finite} derives the finite hook-history rule and its
basic consequences. Section~\ref{sec:cancellation} removes the cutoff
and proves the signed combinatorial rule.
Section~\ref{sec:boundary} evaluates the histories, first in interval
form and then in the single-product form of Theorem~\ref{thm:A}.
Section~\ref{sec:consequences} treats the power-sum expansion, the
one-parameter specialization, stable Grassmannian $K$-homology, and
deformed multipliers. Section~\ref{sec:lowering} deduces the
canonical-basis multiplication rule and the dual lowering rule.
Appendix~\ref{app:tables} records small expansions.

\section{Preliminaries}\label{sec:preliminaries}
This section fixes our conventions for partitions, symmetric functions,
and the Hall-dual canonical bases. We then record the blocking operators
and the finite specialization of the Schur-operator theorem used in the
proof of the main rule.

\subsection{Partitions and skew diagrams}\label{subsec:partitions} 
We begin with conventions for partitions and the elementary facts about skew diagrams. We write a partition as $\lambda=(\lambda_1\geq\lambda_2\geq\cdots\geq0)$, extended by zero parts, and we impose the boundary convention $\lambda_0=\infty$. The set of all partitions is denoted by $\Par$. For a partition $\lambda$, we set the following:
\begin{itemize}
    \item $|\lambda|$ denotes its size,
    \item $\ell(\lambda)$ denotes its length (the number of strictly positive parts),
    \item $\lambda'$ denotes its conjugate.
\end{itemize}
For integers $p\leq q$ we write $[p,q]=\{p,p+1,\dots,q\}$, and
$[n]=[1,n]$. We identify a partition with its English Young diagram, defined by
\[
 [\lambda]=\{(i,j):i\geq1,\ 1\leq j\leq\lambda_i\}.
\]

For $\lambda\subseteq\mu$, the skew diagram $\mu/\lambda=[\mu]\setminus[\lambda]$ has $|\mu/\lambda|=|\mu|-|\lambda|$ boxes. We let $\operatorname{rows}(\mu/\lambda)$ and $\operatorname{cols}(\mu/\lambda)$ denote its number of occupied rows and columns, respectively. 

A box can be added to row $i$ of $\lambda$ if and only if
$\lambda_i<\lambda_{i-1}$; the resulting partition is denoted by
$\lambda+\square_i$. We also put
\begin{equation}\label{eq:Add-def}
 \Add^+(\lambda)=\{i\geq1:0<\lambda_i<\lambda_{i-1}\},
 \qquad
 d(\lambda)=\#\{i\geq1:\lambda_i>\lambda_{i+1}\}.
\end{equation}
Thus $d(\lambda)$ is the number of distinct positive parts of $\lambda$; for
example $d(3,3,2,1)=3$, whereas $\ell(3,3,2,1)=4$.

For every partition $\lambda$, we have
$|\Add^+(\lambda)|=d(\lambda)$. Indeed, the rows with a given positive part $c$ form an interval $[i_c,j_c]$. The
row $i_c$ is the unique row $i$ with $\lambda_i=c<\lambda_{i-1}$, and $j_c$ is
the unique row $j$ with $\lambda_j=c>\lambda_{j+1}$. Conversely, every
$i\in\Add^+(\lambda)$ equals $i_{\lambda_i}$, and every $j$ with
$\lambda_j>\lambda_{j+1}$ has $\lambda_j>0$ and equals $j_{\lambda_j}$. Hence
both sets are in bijection with the set of distinct positive parts.

A skew diagram is \emph{connected} if any two of its boxes are joined by a
chain of boxes of the diagram in which consecutive boxes share an edge. A
\emph{ribbon} is a connected skew diagram containing no $2\times2$ block of
boxes; its height is $\hgt(\mu/\lambda)=\rows(\mu/\lambda)-1$. A
\emph{broken ribbon} is a skew diagram containing no $2\times2$ block,
without any connectivity requirement.

Note that for any $\lambda\subseteq\mu$, the skew shape $\mu/\lambda$ contains no $2\times2$ block if and only if $\mu_{i+1}\leq\lambda_i+1$ for all $i\geq1$. This is immediate: the four boxes $(i,c), (i,c+1), (i+1,c), (i+1,c+1)$ all lie in the skew shape $\mu/\lambda$ if and only if 
\[
 \lambda_i < c < \mu_{i+1}.
\]
Such a $c$ exists if and only if
$\mu_{i+1}\geq\lambda_i+2$.

\subsection{Symmetric functions}\label{subsec:symmetric}

We use the standard notation of \cite[Ch.~I]{Macdonald}, over the coefficient ring $R=\QQ[\alpha,\beta]$. We denote the ring of symmetric functions in the variables $X=(x_1,x_2,\dots)$ over $R$ by $\Lam_R=\bigoplus_{d\geq0}\Lam_R^d$, graded by degree. We let $\widehat{\Lam}_R=\prod_{d\geq0}\Lam_R^d$ denote its degree completion. 

We use the following standard bases and symmetric functions:
\begin{itemize}
    \item $p_k=\sum_ix_i^k$ and $p_\nu=\prod_jp_{\nu_j}$ for the power sums,
    \item $s_\lambda$ and $s_{\lambda/\nu}$ for the Schur and skew Schur functions.
\end{itemize}
We also write $z_\nu=\prod_{i\geq1}i^{m_i(\nu)}m_i(\nu)!$, where $m_i(\nu)$ is the multiplicity of $i$ in $\nu$ \cite{Macdonald}.

The Hall inner product is uniquely determined by $\Hall{s_\lambda}{s_\mu}=\delta_{\lambda\mu}$, or equivalently by $\Hall{p_\lambda}{p_\mu}=\delta_{\lambda\mu}z_\lambda$. Because it orthogonally pairs homogeneous components of equal degree, it naturally extends to a pairing $\widehat{\Lam}_R\times\Lam_R\to R$. Consequently, every completed series $F\in\widehat{\Lam}_R$ satisfies 
\[
 F=\sum_\lambda\Hall{F}{s_\lambda}s_\lambda.
\]

The fundamental involution $\omega$ is the isometric ring automorphism defined by $\omega(p_k)=(-1)^{k-1}p_k$. It satisfies $\omega(s_{\lambda/\nu})=s_{\lambda'/\nu'}$, and it extends degreewise to $\widehat{\Lam}_R$.

The Hopf coproduct $\Delta f=f(X\sqcup Y)$ maps $\widehat{\Lam}_R$ to the completed tensor product $\widehat{\Lam}_R\widehat{\otimes}\widehat{\Lam}_R=\prod_{d,e}\Lam_R^d\otimes\Lam_R^e$. It acts as a ring homomorphism characterized by
\begin{align*}
 \Delta p_k &= p_k\otimes1+1\otimes p_k, \\
 \Delta s_\rho &= \sum_{\nu\subseteq\rho}s_\nu\otimes s_{\rho/\nu},
\end{align*}
as detailed in \cite[Ch.~I, \S5]{Macdonald}. We pair the spaces $\widehat{\Lam}_R\widehat{\otimes}\widehat{\Lam}_R$ and $\Lam_R\otimes\Lam_R$ via the product rule $\Hall{F_1\otimes F_2}{a\otimes b}=\Hall{F_1}{a}\Hall{F_2}{b}$.

We use the following standard facts; see \cite[Ch.~I, \S5]{Macdonald}, and
Example~3 there for (c). For completed elements they follow by applying the
statements for $\Lam_R$ to homogeneous components, since the pairing only involves components of equal degree and is symmetric on each $\Lam^d_R$. 

\begin{lemma}\label{lem:hopf}
\begin{enumerate}[label=\textup{(\alph*)}, itemsep=0.5em]
\item $\Hall{\Delta F}{a\otimes b}=\Hall{F}{ab}$ for all $F\in\widehat{\Lam}_R$ and $a,b\in\Lam_R$.
\item For every $f\in\widehat{\Lam}_R$, there exists a unique $R$-linear adjoint map $f^\perp\colon\Lam_R\to\Lam_R$ such that $\Hall{fH}{q}=\Hall{H}{f^\perp q}$ for all $H\in\widehat{\Lam}_R$ and $q\in\Lam_R$. If $f=\sum_df_d$ with $f_d\in\Lam^d_R$, then $f^\perp q=\sum_d f_d^\perp q$, which evaluates as a finite sum.
\item On $\Lam_R=R[p_1,p_2,\dots]$ we have $p_k^\perp=k\,\partial/\partial p_k$.
This operator extends degreewise to $\widehat{\Lam}_R$, and
$\Hall{p_k^\perp F}{q}=\Hall{F}{p_kq}$ for all $F\in\widehat{\Lam}_R$ and
$q\in\Lam_R$.
\end{enumerate}
\end{lemma}

\emph{Evaluations at scalar alphabets.} For a finite sequence
$B=(b_1,\dots,b_m)$ of elements of $R$, let
$\ev_B,\ev_{-B}\colon\Lam_R\to R$ be the ring homomorphisms with
$\ev_{\pm B}(p_k)=\pm\sum_{i=1}^mb_i^k$. Then $\ev_B(f)=f(b_1,\dots,b_m)$, and
\begin{equation}\label{eq:ev-minus}
 \ev_{-B}(f)=(-1)^d(\omega f)(b_1,\dots,b_m)\quad(f\in\Lam^d_R),
 \qquad
 \ev_{-B}(s_{\rho/\nu})=(-1)^{|\rho/\nu|}s_{\rho'/\nu'}(b_1,\dots,b_m).
\end{equation}
Indeed, $f\mapsto(-1)^{\deg f}(\omega f)(b_1,\dots,b_m)$, extended additively
from homogeneous elements, is a ring homomorphism sending $p_k$ to
$(-1)^k(-1)^{k-1}p_k(b_1,\dots,b_m)=-p_k(b_1,\dots,b_m)$, so it equals
$\ev_{-B}$. The second identity is the first one for $f=s_{\rho/\nu}$, since
$\omega(s_{\rho/\nu})=s_{\rho'/\nu'}$. For an infinite alphabet $B=(b_1,b_2,\ldots)$ of independent
indeterminates, we interpret these evaluations coefficientwise
in the ring of formal sums of monomials of finite support in
the $b_i$. Thus
$\ev_B(p_k)=\sum_{i\geq1}b_i^k$ and
$\ev_{-B}(p_k)=-\sum_{i\geq1}b_i^k$ are well defined, and
\eqref{eq:ev-minus} remains valid. The required constant
specialization is performed only after the tail cancellation
in Proposition~\ref{prop:operator-action}.

\subsection{Canonical and dual canonical Grothendieck functions}
\label{subsec:canonical}

We recall Yeliussizov's two-parameter deformation of the stable
Grothendieck functions and its Hall dual. For parameters $u,v$ and $\ell(\lambda)\leq n$, the finite-variable
canonical stable Grothendieck function is the symmetric formal power
series defined by the following bialternant, expanded at
$x_1=\cdots=x_n=0$ \cite[Def.~3.1]{Yeliussizov}. 
\begin{equation}\label{eq:canonical-determinant}
 G^{(u,v)}_\lambda(x_1,\dots,x_n)
 =\frac{\det\Bigl[x_i^{\lambda_j+n-j}(1+vx_i)^{j-1}(1-ux_i)^{-\lambda_j}\Bigr]_{i,j=1}^n}
       {\prod_{1\leq i<j\leq n}(x_i-x_j)} .
\end{equation}
These series are stable under adjoining a zero variable
\cite[\S3]{Yeliussizov}.  Setting
$G^{(u,v)}_\lambda(x_1,\ldots,x_n)=0$ when $\ell(\lambda)>n$ gives the
stable functions $G^{(u,v)}_\lambda\in\hLam_{\QQ[u,v]}$, and the following hold.
\begin{enumerate}[label=\textup{(G\arabic*)}]
\item\label{G1} $G^{(u,v)}_\lambda=s_\lambda+(\text{terms of degree}>|\lambda|)$,
and every element of $\hLam_{\QQ[u,v]}$ has a unique expansion $\sum_\lambda c_\lambda G^{(u,v)}_\lambda$, convergent degreewise:
since $G^{(u,v)}_\lambda$ has no component of degree $<|\lambda|$, the component of
degree $D$ of the sum involves only the finitely many $\lambda$ with $|\lambda|\leq D$
\cite[\S3, after Def.~3.1]{Yeliussizov}.
\item\label{G2} $G^{(0,0)}_\lambda=s_\lambda$, and $G^{(0,-1)}_\lambda$ is the
stable Grothendieck function of \cite{Buch} \cite[\S3, after
Ex.~3.2]{Yeliussizov}.
\item\label{G3} $\omega\,G^{(u,v)}_\lambda=G^{(v,u)}_{\lambda'}$
\cite[Thm.~5.4]{Yeliussizov}.
\end{enumerate}

For $N\geq0$, let
\[
 \Par_{\leq N}=\{\rho\in\Par:|\rho|\leq N\},
\]
ordered by nondecreasing size, with an arbitrary order among
partitions of the same size. By~\ref{G1}, the finite matrix
\[
 A_N=
 \bigl(\Hall{G^{(-\alpha,-\beta)}_\rho}{s_\nu}\bigr)
 _{\rho,\nu\in\Par_{\leq N}}
\]
is upper unitriangular: its diagonal entries are $1$, and an
off-diagonal entry can be nonzero only when $|\rho|<|\nu|$.
Consequently, $A_N$ is invertible over $R$.

Fix a partition $\lambda$ and take $N=|\lambda|$. Define
\[
 g^{(\alpha,\beta)}_\lambda
 =\sum_{\nu\in\Par_{\leq N}}
   (A_N^{-1})_{\nu,\lambda}\,s_\nu.
\]
Matrix inversion gives the required duality for $|\rho|\leq N$.
For $|\rho|>N$, the pairing vanishes by degree.
Moreover, upper unitriangularity shows that the homogeneous
component of degree $|\lambda|$ is exactly $s_\lambda$.
Thus
\begin{equation}\label{eq:duality}
 \Hall{G^{(-\alpha,-\beta)}_\rho}
       {g^{(\alpha,\beta)}_\lambda}
 =\delta_{\rho\lambda},
 \qquad
 g^{(\alpha,\beta)}_\lambda
 =s_\lambda+(\text{terms of degree}<|\lambda|).
\end{equation}
The same finite matrix argument proves uniqueness among functions
satisfying these conditions.

These are the \emph{dual canonical stable Grothendieck functions}
\cite[Def.~6.1]{Yeliussizov}. They form a basis of $\Lam_R$, and for
$F\in\Lam_R$ we have
$F=\sum_\lambda\Hall{G^{(-\alpha,-\beta)}_\lambda}{F}\,g^{(\alpha,\beta)}_\lambda$,
a finite sum. By~\ref{G2}, $g^{(0,0)}_\lambda=s_\lambda$, and
$g^{(0,1)}_\lambda$ is the dual stable Grothendieck function of Lam and
Pylyavskyy \cite{LamPylyavskyy}, which is Hall dual to $G^{(0,-1)}_\lambda$.

\begin{remark}[Sign conventions]\label{rem:signs}
Several normalizations occur in the literature; see
\cite[Table~1]{IMS25}. In the notation of \cite{Kim-adv, IMS25}, the canonical
functions $G_\lambda(\mathbf x;\balpha,\bbeta)$ with constant parameters
$\alpha_j=\alpha$, $\beta_i=\beta$ are $G^{(-\alpha,-\beta)}_\lambda$ in our
notation, so their dual functions $g_\lambda(\mathbf x;\balpha,\bbeta)$ are our
$g^{(\alpha,\beta)}_\lambda$. For the one-parameter theory we write
\begin{equation}\label{eq:one-parameter-translation}
 \sG^{(\gamma)}_\lambda:=G^{(0,\gamma)}_\lambda,
 \qquad
 \sg^{(\gamma)}_\lambda:=g^{(0,-\gamma)}_\lambda,
 \qquad
 \Hall{\sG^{(\gamma)}_\rho}{\sg^{(\gamma)}_\lambda}=\delta_{\rho\lambda};
\end{equation}
the usual normalization is $\gamma=-1$.
\end{remark}
The involution $\omega$ exchanges the parameters and conjugates
the indexing partition; we record its action on the dual basis.

\begin{lemma}[{\cite[\S6]{Yeliussizov}}]\label{lem:omega-g}
$\omega\,g^{(\alpha,\beta)}_\lambda=g^{(\beta,\alpha)}_{\lambda'}$.
\end{lemma}

\begin{proof}
Since $\omega$ is an isometric involution, \ref{G3} and~\eqref{eq:duality} give
\[
 \Hall{G^{(-\beta,-\alpha)}_\rho}{\omega g^{(\alpha,\beta)}_\lambda}
 =\Hall{\omega G^{(-\beta,-\alpha)}_\rho}{g^{(\alpha,\beta)}_\lambda}
 =\Hall{G^{(-\alpha,-\beta)}_{\rho'}}{g^{(\alpha,\beta)}_\lambda}
 =\delta_{\rho',\lambda}=\delta_{\rho,\lambda'} .
\]
Since also $\omega g^{(\alpha,\beta)}_\lambda=s_{\lambda'}+(\text{lower terms})$,
uniqueness in~\eqref{eq:duality} gives
$\omega g^{(\alpha,\beta)}_\lambda=g^{(\beta,\alpha)}_{\lambda'}$.
\end{proof}

\begin{notation}\label{not:abbrev}
From now on $\alpha,\beta$ are fixed indeterminates, and we abbreviate
\[
 g_\lambda:=g^{(\alpha,\beta)}_\lambda,
 \qquad
 G_\lambda:=G^{(-\alpha,-\beta)}_\lambda,
 \qquad
 b^{(k)}_{\lambda\mu}:=b^{(k)}_{\lambda\mu}(\alpha,\beta),
\]
so that $\Hall{G_\rho}{g_\lambda}=\delta_{\rho\lambda}$. For a symmetric
function $F\in\Lam_R$ we write $[g_\mu]F$ for the coefficient of $g_\mu$ in
$F$.
\end{notation}

Both bases are homogeneous, but for opposite gradings of the parameters.

\begin{lemma}\label{lem:grading}
Write $G_\rho=\sum_\nu a_{\rho\nu}s_\nu$ and $g_\lambda=\sum_\nu c_{\lambda\nu}s_\nu$
with $a_{\rho\nu},c_{\lambda\nu}\in R$.
\begin{enumerate}[label=\textup{(\alph*)}]
\item $a_{\rho\nu}$ is a homogeneous polynomial of degree $|\nu|-|\rho|$.
\item $c_{\lambda\nu}$ is a homogeneous polynomial of degree $|\lambda|-|\nu|$.
\end{enumerate}
In other words, $G_\rho$ is homogeneous of degree $|\rho|$ if
$\deg x_i=1$ and $\deg\alpha=\deg\beta=-1$, and $g_\lambda$ is homogeneous
of degree $|\lambda|$ if $\deg x_i=\deg\alpha=\deg\beta=1$.
\end{lemma}

\begin{proof}
Let $t$ be an indeterminate. A polynomial $f\in R$ is homogeneous of
degree $m$ if and only if $f(\alpha/t,\beta/t)=t^{-m}f(\alpha,\beta)$ in
$R[t,t^{-1}]$.

(a) In~\eqref{eq:canonical-determinant} with $(u,v)=(-\alpha,-\beta)$,
substitute $x_i\mapsto tx_i$, $\alpha\mapsto\alpha/t$ and
$\beta\mapsto\beta/t$. The factors $(1-\beta x_i)$ and $(1+\alpha x_i)$ are
unchanged, the $j$-th column of the matrix is multiplied by
$t^{\rho_j+n-j}$, and the denominator by $t^{\binom n2}$. Hence
$G_\rho(tx;\alpha/t,\beta/t)=t^{|\rho|}G_\rho(x;\alpha,\beta)$ for every $n$,
and therefore in $\hLam_{R[t,t^{-1}]}$. Since $s_\nu(tx)=t^{|\nu|}s_\nu(x)$,
comparing coefficients of $s_\nu$ gives
$a_{\rho\nu}(\alpha/t,\beta/t)=t^{|\rho|-|\nu|}a_{\rho\nu}(\alpha,\beta)$.

(b) Since the Schur functions are orthonormal, the duality~\eqref{eq:duality}
reads
\begin{equation}\label{eq:duality-coeff}
 \sum_\nu a_{\rho\nu}\,c_{\lambda\nu}=\delta_{\rho\lambda}
 \qquad\text{for all }\rho,\lambda,
\end{equation}
a finite sum, and $c_{\lambda\lambda}=1$, $c_{\lambda\nu}=0$ unless $\nu=\lambda$
or $|\nu|<|\lambda|$. Apply the ring homomorphism $R\to R[t,t^{-1}]$,
$\alpha\mapsto\alpha/t$, $\beta\mapsto\beta/t$, to~\eqref{eq:duality-coeff}
and use (a). With
$c'_{\lambda\nu}:=t^{|\lambda|-|\nu|}\,c_{\lambda\nu}(\alpha/t,\beta/t)$ this gives
\[
 \sum_\nu t^{|\rho|-|\nu|}a_{\rho\nu}\cdot t^{|\nu|-|\lambda|}c'_{\lambda\nu}
 =\delta_{\rho\lambda},
 \qquad\text{that is,}\qquad
 \sum_\nu a_{\rho\nu}\,c'_{\lambda\nu}=\delta_{\rho\lambda},
\]
because $t^{|\rho|-|\lambda|}=1$ when $\rho=\lambda$. Moreover
$c'_{\lambda\lambda}=1$, and $c'_{\lambda\nu}=0$ unless $\nu=\lambda$ or
$|\nu|<|\lambda|$. The solution of~\eqref{eq:duality-coeff} with these
properties is unique, since the matrix $(a_{\rho\nu})$ is unitriangular with
respect to size by~\ref{G1}; this is the uniqueness in~\eqref{eq:duality},
and its proof works over any coefficient ring. Hence
$c'_{\lambda\nu}=c_{\lambda\nu}$, that is,
$c_{\lambda\nu}(\alpha/t,\beta/t)=t^{|\nu|-|\lambda|}c_{\lambda\nu}(\alpha,\beta)$.
\end{proof}

\emph{Coproduct-skew functions.} Using~\ref{G1} componentwise in the second tensor factor of $\Delta G_{\mu}$, there are unique elements
$G_{\mu\sk\lambda}\in\hLam_R$ such that
\begin{equation}
\label{eq:coproduct-skew}
 \Delta G_\mu=\sum_{\lambda}G_{\mu\sk\lambda}\otimes G_\lambda .
\end{equation}

These are the skew functions defined by Schur operators in
\cite[\S5.4, Eq.(34)]{Yeliussizov} (for $(\alpha,\beta)=(0,1)$ they are the
functions $G_{\mu/\!\!/\lambda}$ of Buch \cite{Buch}), and they are studied in
\cite[\S4.1]{IMS24}. The double slash distinguishes them from the tableau
skew functions $G_{\mu/\lambda}$, which differ from them in general.

These coproduct-skew functions encode multiplication in the
dual basis through the Hall pairing.

\begin{lemma}\label{lem:skew-pairing}
For every $f\in\Lam_R$ and every partition $\lambda$,
\[
 f\,g_\lambda=\sum_\mu\Hall{G_{\mu\sk\lambda}}{f}\,g_\mu ,
\]
a finite sum. In particular $b^{(k)}_{\lambda\mu}=\Hall{G_{\mu\sk\lambda}}{p_k}$.
\end{lemma}

\begin{proof}
Since $fg_\lambda\in\Lam_R$, we have
$fg_\lambda=\sum_\mu\Hall{G_\mu}{fg_\lambda}g_\mu$ with finitely many nonzero
terms. By Lemma~\ref{lem:hopf}(a) and~\eqref{eq:coproduct-skew},
\[
 \Hall{G_\mu}{fg_\lambda}=\Hall{\Delta G_\mu}{f\otimes g_\lambda}
 =\sum_\nu\Hall{G_{\mu\sk\nu}}{f}\Hall{G_\nu}{g_\lambda}
 =\Hall{G_{\mu\sk\lambda}}{f}.
\]
The middle sum may be computed termwise: only bidegrees $(d_1,d_2)$ with
$d_1\leq\deg f$ and $d_2\leq|\lambda|$ contribute to the pairing, and
by~\ref{G1} the function $G_\nu$ has no component of degree $<|\nu|$, so in
each such bidegree only the finitely many $\nu$ with $|\nu|\leq d_2$
contribute.
\end{proof}

\begin{remark}\label{rem:hopf-skew}
The identity $b^{(k)}_{\lambda\mu}=\Hall{G_{\mu\sk\lambda}}{p_k}$ is a
cutoff-free expression, but expanding $p_k$ by~\eqref{eq:hook-expansion}
hides an alternating cancellation. Theorems~\ref{thm:connected-rule}
and~\ref{thm:boundary} carry out this cancellation explicitly.
\end{remark}

\subsection{Blocking operators}\label{subsec:operators}

The operators of Iwao, Motegi and Scrimshaw act on the free module spanned
by partitions. We describe them for constant parameters, together with
the two elementary properties used below.

Let $R[\Par]$ be the free $R$-module with basis $\Par$, and for
$F\in R[\Par]$ let $[\mu]F$ be the coefficient of $\mu$. For $i\geq1$ and a
partition $\nu$ define the \emph{stay scalar}
\begin{equation}\label{eq:stay-scalar}
 c_i(\nu)=
 \begin{cases}
  -\alpha,&0<\nu_i<\nu_{i-1},\\
  0,&0=\nu_i<\nu_{i-1},\\
  \beta,&\nu_i=\nu_{i-1},
 \end{cases}
\end{equation}
and define $U_i\colon R[\Par]\to R[\Par]$ by
\begin{equation}\label{eq:Ui}
 U_i\nu=\ind_{\{\nu_i<\nu_{i-1}\}}\,(\nu+\square_i)+c_i(\nu)\,\nu .
\end{equation}
Since $\nu_0=\infty$, $U_1\nu$ always has an addition term. If
$\nu_1>0$, it also has the stay term $-\alpha\nu$. A \emph{move}
of $U_i$ at $\nu$ is either the \emph{addition} $\nu\mapsto\nu+\square_i$,
which is admissible if $\nu_i<\nu_{i-1}$ and has weight $1$, or the
\emph{stay} $\nu\mapsto\nu$, which is admissible if $c_i(\nu)\neq0$ and has
weight $c_i(\nu)$. Thus $U_i\nu$ is the sum over the admissible moves of
their weights times their results. These operators are the constant
specializations of the blocking operators of \cite[\S3.1]{IMS25}; for
$\alpha=0$ they were introduced in \cite[\S6]{Iwao20}.

\begin{lemma}\label{lem:local}
The admissibility and the weight of a move of $U_i$ at $\nu$ depend only on
$(\nu_{i-1},\nu_i)$, and a move of $U_i$ changes only $\nu_i$. In
particular $U_iU_j=U_jU_i$ for $|i-j|\geq2$.
\end{lemma}

\begin{proof} 
The first statement is immediate from~\eqref{eq:stay-scalar}
and~\eqref{eq:Ui}. Using~\eqref{eq:Ui} we obtain, for every partition $\nu$ (terms with a vanishing
indicator $\ind_{\{\nu_i<\nu_{i-1}\}}$ are omitted, so $\nu+\square_i$ only occurs when it exists),
\begin{align*}
 U_jU_i\nu&=U_j\{\ind_{\{\nu_i<\nu_{i-1}\}}(\nu+\square_i)+c_i(\nu)\nu\}\\
 &= \ind_{\{\nu_i<\nu_{i-1}\}}[\ind_{\{(\nu+\square_i)_j<(\nu+\square_i)_{j-1}\}} (\nu + \square_i+\square_j) + c_j(\nu+\square_i)(\nu+\square_i) ]\\& + c_i(\nu)[ \ind_{\{\nu_j<\nu_{j-1}\}}\,(\nu+\square_j)+c_j(\nu)\,\nu ].
 \end{align*}
Suppose $\nu_i < \nu_{i-1}$, then $ \nu+\square_i$ exists. In this case, since $|i-j| \geq 2$, $(\nu+\square_i)_j=\nu_j$ and $(\nu+\square_i)_{j-1}=\nu_{j-1}$. Thus we have
\begin{align*}
U_jU_i(\nu)&=\ind_{\{\nu_i<\nu_{i-1}\}}\ind_{\{\nu_j<\nu_{j-1}\}} (\nu + \square_i+\square_j)+
\ind_{\{\nu_i<\nu_{i-1}\}}c_j(\nu)\,(\nu+\square_i)\\ &+\ind_{\{\nu_j<\nu_{j-1}\}}c_i(\nu)\,(\nu+\square_j)+c_i(\nu)c_j(\nu)\nu.
\end{align*}
Interchanging $i$ and $j$, we conclude that $U_iU_j=U_jU_i$.
\end{proof}
Since each move either adds one box or leaves the partition
unchanged, every word acts triangularly with respect to
containment.

\begin{lemma}\label{lem:triangular}
For a word $w=w_1\cdots w_N$, put $U_w=U_{w_1}\cdots U_{w_N}$. For every
partition $\nu$, $[\mu]\,U_w\nu\neq0$ implies $\nu\subseteq\mu$ and
$|\mu/\nu|\leq N$, and
$[\nu]\,U_w\nu=\prod_{j=1}^Nc_{w_j}(\nu)$.
\end{lemma}

\begin{proof}
By~\eqref{eq:Ui}, each $U_i$ maps a partition $\kappa$ to $c_i(\kappa)\kappa$
plus a multiple of a partition containing $\kappa$ with one more box.
Iterating $N$ times, every partition occurring in $U_w\nu$ contains $\nu$
and has at most $N$ more boxes. A term returns to $\nu$ only if every step
is a stay; then all intermediate partitions equal $\nu$, and the
coefficient is $\prod_jc_{w_j}(\nu)$.
\end{proof}

\emph{Noncommutative Schur functions.} Fix $n\geq1$ and write
$\bfU_n=(U_1,\dots,U_n)$. The operators satisfy (see 
\cite[Lemma~3.1]{IMS25})
\begin{equation}
    \label{eq:Knuth-1}
    U_iU_j=U_jU_i\quad(|i-j|\geq2),
\end{equation}
\begin{equation}
    \label{eq:Knuth-2}
    (U_i+U_{i+1})U_{i+1}U_i=U_{i+1}U_i(U_i+U_{i+1}).
\end{equation}
The first relation~\eqref{eq:Knuth-1} is Lemma~\ref{lem:local}, which is
the non-local commutativity relation \cite[\S1, (1.1)]{FominGreene}, and the
second~\eqref{eq:Knuth-2} is the local commutativity relation \cite[\S1, (1.3)]{FominGreene}. Using \eqref{eq:Knuth-1}, for $ i<j<k$, we have
\begin{equation}
    \label{eq:commute-1}
    U_iU_kU_j=U_kU_iU_j \text{ and }
    U_jU_iU_k=U_jU_kU_i,
\end{equation}
and for $i\geq 1,k\geq 2$, we obtain
\begin{equation}
    \label{eq:commute-2}
    U_{i+k}U_i(U_i+U_{i+k})=(U_i+U_{i+k})U_{i+k}U_i.
\end{equation}
In view of \eqref{eq:Knuth-2}, \eqref{eq:commute-1} and \eqref{eq:commute-2}, \cite[Lemma~3.1]{FominGreene} implies that the noncommutative
elementary symmetric functions
$e_m(\bfU_n)=\sum_{n\geq i_1>\cdots>i_m\geq1}U_{i_1}\cdots U_{i_m}$ commute
pairwise. Consequently, there is a unique ring homomorphism
\[
 \Phi_n\colon\Lam_R\longrightarrow\operatorname{End}_R(R[\Par]),
 \qquad \Phi_n(e_m)=e_m(\bfU_n),
\]
whose image is commutative. For every partition $\theta$, we write
$s_\theta(\bfU_n):=\Phi_n(s_\theta)$. The dual Jacobi--Trudi identity is
\[
 s_\theta
 =\det\bigl(e_{\theta'_i-i+j}\bigr)_{1\leq i,j\leq m},
 \qquad m\geq\ell(\theta')=\theta_1,
\]
where $e_0=1$ and $e_j=0$ for $j<0$
\cite[Ch.~I, \S3]{Macdonald}. Applying the ring homomorphism
$\Phi_n$ therefore gives
\[
 s_\theta(\bfU_n)
 =\det\bigl(e_{\theta'_i-i+j}(\bfU_n)\bigr)_{1\leq i,j\leq m}.
\]
According to \cite[Lemma~3.2]{FominGreene}, we have the tableau formula
\begin{equation}\label{eq:FGtableau}
 s_\theta(\bfU_n)=\sum_{T\in\SSYT_n(\theta)}U_{w(T)} ,
\end{equation}
where $\SSYT_n(\theta)$ is the set of semistandard Young tableaux of shape
$\theta$ with entries in $[n]=\{1,\dots,n\}$, and $w(T)$ is the \emph{column
reading word}, obtained by reading each column from bottom to top and the
columns from left to right. Recall that $U_{w(T)}=U_{w_1}\cdots U_{w_N}$, so the
last letter of $w(T)$ acts first.

\emph{Supersymmetric evaluation.} For $n\geq1$ and $m\geq0$, let
$\ev_{-\bbeta_m}\colon\Lam_R\to R$ be the evaluation at the negative of the
scalar alphabet $\bbeta_m=(\beta,\dots,\beta)$ ($m$ entries), and define the
ring homomorphism
\begin{equation}\label{eq:Phi-super}
 \Phi_{n,m}=(\Phi_n\otimes\ev_{-\bbeta_m})\circ\Delta\colon
 \Lam_R\to\operatorname{End}_R(R[\Par]).
\end{equation}
By $\Delta s_\rho=\sum_\nu s_\nu\otimes s_{\rho/\nu}$ and~\eqref{eq:ev-minus},
\begin{equation}\label{eq:super-schur}
 s_\rho(\bfU_n/\bbeta_m):=\Phi_{n,m}(s_\rho)
 =\sum_{\nu\subseteq\rho}(-1)^{|\rho/\nu|}\,s_{\rho'/\nu'}(\bbeta_m)\,s_\nu(\bfU_n),
\end{equation}
and since $\Delta p_k=p_k\otimes1+1\otimes p_k$,
\begin{equation}\label{eq:superpower}
 \Phi_{n,m}(p_k)=p_k(\bfU_n)-m\beta^k,
 \qquad
 p_k(\bfU_n):=\Phi_n(p_k)=\sum_{a=0}^{k-1}(-1)^as_{(k-a,1^a)}(\bfU_n).
\end{equation}

\subsection{The Schur-operator theorem}\label{subsec:IMS}

The following proposition is the starting point of the paper. It is the
constant-parameter case of a theorem of Iwao, Motegi and Scrimshaw \cite{IMS25}; since
that theorem is stated with infinitely many parameters $\beta_i$, whose
constant specialization is not defined termwise, we include the truncation
argument that makes the specialization legitimate.

\begin{proposition}\label{prop:operator-action}
Let $\lambda$ be a partition, $N\geq0$, and $r\geq\ell(\lambda)+N$. For every
$f\in\Lam^N_R$ and every partition $\mu$,
\[
 [g_\mu]\bigl(f\,g_\lambda\bigr)=[\mu]\,\Phi_{r+1,r}(f)\,\lambda .
\]
\end{proposition}

\begin{proof}
Both sides are $R$-linear in $f$, so it suffices to take $f=s_\rho$ with
$|\rho|=N$. We use \cite{IMS25} with refined parameters: indeterminates
$\balpha=(\alpha_1,\alpha_2,\dots)$ with $\alpha_0:=0$, and
$\bbeta=(\beta_1,\beta_2,\dots)$. The refined blocking operator
$U^{\natural}_i$ acts on $\nu$ by the addition (if $\nu_i<\nu_{i-1}$) plus
the stay with weight $-\alpha_{\nu_i}$ if $\nu_i<\nu_{i-1}$, and with weight
$\beta_{i-1}$ if $\nu_i=\nu_{i-1}$. Specializing $\alpha_j=\alpha$ $(j\geq1)$ and
$\beta_i=\beta$ gives~\eqref{eq:Ui}; note that the weight $-\alpha_0=0$
reproduces the second case of~\eqref{eq:stay-scalar}.

Let
$G_{\mu\sk\lambda}(\mathbf x;\balpha,\bbeta)$ and the vectors
${}^{[\balpha,\bbeta]}\langle\lambda|$ and $|\mu\rangle^{[\balpha,\bbeta]}$ be as
in \cite[\S2.3]{IMS25}, and let $a_k$ ($k\geq1$) be the current operators.
By \cite[(2.6), (2.7)]{IMS25},
\[
 {}^{[\balpha,\bbeta]}\langle\lambda|\mu\rangle^{[\balpha,\bbeta]}=\delta_{\lambda\mu},
 \qquad
 G_{\mu\sk\lambda}(\mathbf x;\balpha,\bbeta)
 ={}^{[\balpha,\bbeta]}\langle\lambda|\,e^{H(\mathbf x)}\,|\mu\rangle^{[\balpha,\bbeta]},
\]
where $H(\mathbf x)=\sum_{k\geq1}p_k(\mathbf x)a_k/k$. For each partition
$\rho$, let $S_\rho\in\QQ[p_1,p_2,\dots]$ denote the unique polynomial that
expresses the Schur function in terms of the power sums, so that
$s_\rho(\mathbf y)=S_\rho(p_1(\mathbf y),p_2(\mathbf y),\dots)$. We have
\begin{align*}
 \exp\Bigl(\sum_{k\geq1}p_k(\mathbf x)p_k(\mathbf y)/k\Bigr)
 &=\prod_{k\geq1}\exp\bigl(p_k(\mathbf x)p_k(\mathbf y)/k\bigr)\\
 &=\prod_{k\geq1}\sum_{m_k=0}^{\infty}
   \frac{\bigl(p_k(\mathbf x)p_k(\mathbf y)\bigr)^{m_k}}{k^{m_k}m_k!}\\
 &=\sum_{\lambda}z_{\lambda}^{-1}p_{\lambda}(\mathbf{x})p_{\lambda}(\mathbf{y}).
\end{align*}
Equating the power-sum and Schur expansions of the Cauchy kernel
\cite[Ch.~I, (4.1) and (4.3)]{Macdonald} with the above identity yields
\[
 \exp\Bigl(\sum_{k\geq1}p_k(\mathbf x)p_k/k\Bigr)
 =\sum_\rho s_\rho(\mathbf x)S_\rho(p_1,p_2,\dots).
\]
Since the operators $a_k$ ($k\geq1$) commute pairwise, substituting
$p_k\mapsto a_k$ gives
$e^{H(\mathbf x)}=\sum_\rho s_\rho(\mathbf x)S_\rho(a_1,a_2,\dots)$. The matrix
elements
$C_\rho:={}^{[\balpha,\bbeta]}\langle\lambda|S_\rho(a_1,a_2,\dots)|\mu\rangle^{[\balpha,\bbeta]}$
do not depend on $\mathbf x$, and matrix elements are linear, so
\[
 G_{\mu\sk\lambda}(\mathbf x;\balpha,\bbeta)=\sum_\rho C_\rho\,s_\rho(\mathbf x),
\]
where in each degree $D$ in $\mathbf x$ only the finitely many $\rho$ with
$|\rho|=D$ occur. Since the Hall inner product is taken in $\mathbf x$ and
$\{s_\rho(\mathbf x)\}$ is orthonormal ($\Hall{s_\sigma}{s_\rho}=\delta_{\sigma\rho}$),
pairing with $s_\rho$ gives
\[
 \Hall{G_{\mu\sk\lambda}(\mathbf x;\balpha,\bbeta)}{s_\rho}
 =C_\rho
 ={}^{[\balpha,\bbeta]}\langle\lambda|\,S_\rho(a_1,a_2,\dots)\,|\mu\rangle^{[\balpha,\bbeta]}.
\]
By \cite[Thm.~3.2]{IMS25},
${}^{[\balpha,\bbeta]}\langle\lambda|S_\rho(a_1,a_2,\dots)
={}^{[\balpha,\bbeta]}\langle s_\rho(\bfU^\natural/\bbeta)\lambda|$, where
\begin{equation}\label{eq:super}
 s_\rho(\bfU^\natural/\bbeta)\lambda
 =\sum_{\nu\subseteq\rho}(-1)^{|\rho/\nu|}s_{\rho'/\nu'}(\beta_1,\beta_2,\dots)\,
 s_\nu(U^\natural_1,U^\natural_2,\dots)\lambda
\end{equation}
is formed with infinite alphabets. Since each $U^\natural_i$ adds at most one
box, $s_\rho(\bfU^\natural/\bbeta)\lambda$ is a finite linear combination of
partitions $\eta\supseteq\lambda$ with $|\eta/\lambda|\leq N$, whose
coefficients are formal power series in the $\beta_i$. By the orthonormality
${}^{[\balpha,\bbeta]}\langle\eta|\mu\rangle^{[\balpha,\bbeta]}=\delta_{\eta\mu}$,
\begin{equation}\label{eq:refined-identity}
 \Hall{G_{\mu\sk\lambda}(\mathbf x;\balpha,\bbeta)}{s_\rho}
 =[\mu]\,s_\rho(\bfU^\natural/\bbeta)\,\lambda .
\end{equation}

Every partition obtained from $\lambda$ by at most $N$ moves has at most
$\ell(\lambda)+N\leq r$ nonzero parts. For such a partition $\nu$ and every
$i\geq r+2$ we have $\nu_{i-1}=\nu_i=0$, hence
$U^\natural_i\nu=\beta_{i-1}\nu$. Let $\theta$ be a partition with
$|\theta|\leq N$, and let $T$ be a semistandard tableau of shape $\theta$ with
positive entries. The entries $\leq r+1$ of $T$ form a subtableau
$T_{\leq}$ of some shape $\kappa\subseteq\theta$, and the remaining entries
form a skew tableau $T_>$ of shape $\theta/\kappa$ with entries $\geq r+2$.
In every column the entries $\leq r+1$ lie above the others, so deleting
the letters $\geq r+2$ from $w(T)$ leaves $w(T_{\leq})$. When
$U^\natural_{w(T)}$ is applied to $\lambda$, every intermediate partition is
obtained from $\lambda$ by fewer than $N$ moves, so each letter $j\geq r+2$
acts as the scalar $\beta_{j-1}$. Therefore
\[
 U^\natural_{w(T)}\lambda=\Bigl(\prod_{j\in T_>}\beta_{j-1}\Bigr)\,
 U^\natural_{w(T_{\leq})}\lambda .
\]
Summing over $T$ and using the tableau formula~\eqref{eq:FGtableau} for the
alphabets $U^\natural_1,U^\natural_2,\dots$ and
$U^\natural_1,\dots,U^\natural_{r+1}$,
\[
 s_\theta(U^\natural_1,U^\natural_2,\dots)\lambda
 =\sum_{\kappa\subseteq\theta}s_{\theta/\kappa}(\beta_{r+1},\beta_{r+2},\dots)\,
  s_\kappa(U^\natural_1,\dots,U^\natural_{r+1})\lambda .
\]
Write $\bbeta_{>r}=(\beta_{r+1},\beta_{r+2},\dots)$ and
$\bbeta_{\leq r}=(\beta_1,\dots,\beta_r)$, and let
$\epsilon=(\ev_{\bbeta_{>r}}\otimes\ev_{-\bbeta})\circ\Delta$ (evaluations at the
infinite alphabets, as in Section~\ref{subsec:symmetric}). It is a ring
homomorphism with
$\epsilon(p_k)=\sum_{i>r}\beta_i^k-\sum_{i\geq1}\beta_i^k=-\sum_{i\leq r}\beta_i^k$,
that is, $\epsilon=\ev_{-\bbeta_{\leq r}}$. Substituting the previous identity into
\eqref{eq:super} and collecting the coefficient of
$s_\kappa(U^\natural_1,\dots,U^\natural_{r+1})\lambda$ in
$s_\rho(\bfU^\natural/\bbeta)\lambda$ gives
\[
 \sum_{\kappa\subseteq\nu\subseteq\rho}(-1)^{|\rho/\nu|}
 s_{\rho'/\nu'}(\bbeta)\,s_{\nu/\kappa}(\bbeta_{>r})
 =\sum_{\kappa\subseteq\nu\subseteq\rho}\ev_{\bbeta_{>r}}(s_{\nu/\kappa})\,\ev_{-\bbeta}(s_{\rho/\nu})
 =\epsilon(s_{\rho/\kappa}),
\]
where the first equality is \eqref{eq:ev-minus}, and the second follows from
\[
 \Delta s_{\rho/\kappa}=\sum_{\kappa\subseteq\nu\subseteq\rho}s_{\nu/\kappa}\otimes s_{\rho/\nu}.
\]
Indeed, for partitions $\sigma,\tau$, with $c^\nu_{\kappa\sigma}=\Hall{s_{\nu/\kappa}}{s_\sigma}$,
\begin{align*}
 \Hall{\Delta s_{\rho/\kappa}}{s_\sigma\otimes s_\tau}
 &=\Hall{s_{\rho/\kappa}}{s_\sigma s_\tau}=\Hall{s_\rho}{s_\kappa s_\sigma s_\tau}
 =\sum_\nu c^\nu_{\kappa\sigma}\Hall{s_\rho}{s_\nu s_\tau}
 =\sum_\nu c^\nu_{\kappa\sigma}\Hall{s_{\rho/\nu}}{s_\tau}\\
 &=\sum_{\kappa\subseteq\nu\subseteq\rho}\Hall{s_{\nu/\kappa}}{s_\sigma}\Hall{s_{\rho/\nu}}{s_\tau}.
\end{align*}
By \eqref{eq:ev-minus} again,
$\epsilon(s_{\rho/\kappa})=(-1)^{|\rho/\kappa|}s_{\rho'/\kappa'}(\beta_1,\dots,\beta_r)$. Thus
\[
 s_\rho(\bfU^\natural/\bbeta)\lambda
 =\sum_{\kappa\subseteq\rho}(-1)^{|\rho/\kappa|}
 s_{\rho'/\kappa'}(\beta_1,\dots,\beta_r)\,
 s_\kappa(U^\natural_1,\dots,U^\natural_{r+1})\lambda.
\]
The right-hand side involves only finitely many parameters and is
polynomial in them.

Now we specialize $\alpha_j=\alpha$ for $j\geq1$ and
$\beta_i=\beta$ for $i\geq1$. Then the right-hand side of~\eqref{eq:refined-identity}
becomes $[\mu]\,s_\rho(\bfU_{r+1}/\bbeta_r)\lambda=[\mu]\,\Phi_{r+1,r}(s_\rho)\lambda$,
by~\eqref{eq:super-schur}. On the left-hand side, the functions
$G_\mu(\mathbf x;\balpha,\bbeta)$ specialize to $G_\mu=G^{(-\alpha,-\beta)}_\mu$
(Remark~\ref{rem:signs}). Their coproduct-skew functions are characterized
by the branching rule
$G_\mu(\mathbf x,\mathbf y)=\sum_\nu G_{\mu\sk\nu}(\mathbf y)G_\nu(\mathbf x)$
\cite[Prop.~2.3]{IMS25}, that is, by~\eqref{eq:coproduct-skew}, so they
specialize to our $G_{\mu\sk\lambda}$. The left-hand side therefore becomes
$\Hall{G_{\mu\sk\lambda}}{s_\rho}$, which is $[g_\mu](s_\rho g_\lambda)$ by
Lemma~\ref{lem:skew-pairing}.
\end{proof}

\section{The finite hook-history rule}\label{sec:finite}
In this section we first turn the Schur-operator theorem into a finite coefficient
formula. Its terms are histories of semistandard hook words; their
additions determine the final partition and their stays record the
parameter weights.

\subsection{Hook words}

We first describe the tableaux of hook shape and the order in which the
letters of their reading words act.

For $k\geq1$ and $0\leq a<k$ put $H_{k,a}=(k-a,1^a)$ and $t=k-a$.

\begin{lemma}\label{lem:hook-words}
For every $n\geq1$, the map $T\mapsto w(T)$ is a bijection from
$\SSYT_n(H_{k,a})$ onto the set of words
\begin{equation}\label{eq:hook-word}
 w=d_a\cdots d_1\,c_1c_2\cdots c_t,
 \qquad
 n\geq d_a>\cdots>d_1>c_1\leq c_2\leq\cdots\leq c_t\leq n .
\end{equation}
\end{lemma}

\begin{proof}
A tableau $T$ of shape $H_{k,a}$ has first column $c_1<d_1<\cdots<d_a$, read
from top to bottom, and first row $c_1\leq c_2\leq\cdots\leq c_t$. Its column
reading word is the first column read upwards, $d_a\cdots d_1c_1$, followed
by the one-box columns $c_2,\dots,c_t$. Conversely, every
word~\eqref{eq:hook-word} arises from exactly one such tableau.
\end{proof}

A word of the form~\eqref{eq:hook-word}, with arbitrary positive letters, is
called a \emph{hook word} of length $k$. We write $\leg(w)=a$, call
$\Dec(w)=\{d_1,\dots,d_a\}$ its decreasing part and
$\Inc(w)=(c_1,\dots,c_t)$ its increasing part (a multiset), and let
$\supp(w)$ be its set of letters. Note that $c_1=\min\supp(w)$. Since the
last letter acts first, the operators of $U_w$ act on a vector in the
\emph{chronological order}
\begin{equation}\label{eq:two-pass}
 \underbrace{c_t,\ c_{t-1},\ \dots,\ c_1}_{\text{first pass (weakly decreasing)}}
 ,\ \underbrace{d_1,\ d_2,\ \dots,\ d_a}_{\text{second pass (strictly increasing)}} .
\end{equation}

\subsection{Histories}

Applying a word in the operators $U_i$ to a partition and expanding
by~\eqref{eq:Ui} amounts to choosing a move for each letter. We record such
choices as follows.

\begin{definition}\label{def:history}
Let $w=w_1\cdots w_N$ be a word and $\lambda$ a partition. A \emph{history}
of $w$ from $\lambda$ is a choice of move $h_j$ of $U_{w_j}$ for each
position $j\in\{1,\ldots,N\}$. Applying the moves in chronological order
\[
 \lambda = \nu^{(0)} \xrightarrow{h_N} \nu^{(1)} \xrightarrow{h_{N-1}} \cdots \xrightarrow{h_1} \nu^{(N)},
\]
each move $h_j$ must be admissible at the partition $\nu^{(N-j)}$ to which
it is applied. The final partition $\nu^{(N)}$ is denoted $\enD(h)$.
\end{definition}

Let $\#_\alpha(h)$ and $\#_\beta(h)$ be the numbers of stays of
weight $-\alpha$ and $\beta$, respectively, and put
\[
 \wt(h)=(-\alpha)^{\#_\alpha(h)}\beta^{\#_\beta(h)} .
\]
\begin{example}\label{ex:histories}
Let $w=312$ and $\lambda=(3,2,0)$, so $U_w=U_3U_1U_2$. The operators act in
the order $U_2,U_1,U_3$ (the rightmost acts first), and the moves are applied as in
\[
 \lambda=\nu^{(0)}\xrightarrow{h_3}\nu^{(1)}\xrightarrow{h_2}\nu^{(2)}\xrightarrow{h_1}\nu^{(3)},
\]
where $h_3,h_2,h_1$ are moves of $U_2,U_1,U_3$, respectively. At each step at most
one stay is admissible, namely the one of weight $c_i(\nu)$ when $c_i(\nu)\neq0$.
\begin{itemize}
\item $U_2$ acts at $(3,2,0)$, where $0<\nu_2<\nu_1$: it can add or stay (weight $-\alpha$).
\item $U_1$ acts next, where $0<\nu_1<\nu_0=\infty$: it can add or stay (weight $-\alpha$).
\item $U_3$ acts last, where $0=\nu_3<\nu_2$: it can only add, since $c_3(\nu)=0$.
\end{itemize}
Hence $w$ has $2\cdot2\cdot1=4$ histories from $\lambda$:
\[
\renewcommand{\arraystretch}{1.2}
\begin{array}{lll}
(h_3,h_2,h_1) & \nu^{(0)}\to\nu^{(1)}\to\nu^{(2)}\to\nu^{(3)} & \wt(h)\\ \hline
\text{add, add, add} & (3,2,0)\to(3,3,0)\to(4,3,0)\to(4,3,1) & 1\\
\text{add, stay, add} & (3,2,0)\to(3,3,0)\to(3,3,0)\to(3,3,1) & -\alpha\\
\text{stay, add, add} & (3,2,0)\to(3,2,0)\to(4,2,0)\to(4,2,1) & -\alpha\\
\text{stay, stay, add} & (3,2,0)\to(3,2,0)\to(3,2,0)\to(3,2,1) & \alpha^2
\end{array}
\]
Therefore
\[
U_3U_1U_2\,(3,2,0)=(4,3,1)-\alpha\,(3,3,1)-\alpha\,(4,2,1)+\alpha^2\,(3,2,1).
\]
\end{example}

The coefficients of $U_w\lambda$ are sums over histories.

\begin{lemma}\label{lem:histories}
For every word $w$ of length $N$ and every partition $\lambda$,
$[\mu]\,U_w\lambda=\sum_h\wt(h)$, the sum over histories of $w$ from $\lambda$
with $\enD(h)=\mu$. Every history satisfies $\enD(h)\supseteq\lambda$, and the
number of its additions is $|\enD(h)/\lambda|$, so
$\#_\alpha(h)+\#_\beta(h)=N-|\enD(h)/\lambda|$.
\end{lemma}

\begin{proof}
Expanding $U_{w_1}\cdots U_{w_N}\lambda$ with~\eqref{eq:Ui}, each term
corresponds to a choice of an admissible move for each letter, applied in
the order $w_N,\dots,w_1$, and its coefficient is the product of the
weights. Each addition adds one box and each stay adds none, which gives
the remaining statements.
\end{proof}

For $r\geq0$ define
\begin{equation}\label{eq:Mdef}
 M^r_{k,a}(\lambda,\mu)
 =\sum_{T\in\SSYT_{r+1}(H_{k,a})}\
 \sum_{\substack{h\text{ history of }w(T)\text{ from }\lambda\\ \enD(h)=\mu}}\wt(h).
\end{equation}
By~\eqref{eq:FGtableau} and Lemma~\ref{lem:histories},
$M^r_{k,a}(\lambda,\mu)=[\mu]\,s_{H_{k,a}}(\bfU_{r+1})\lambda$.

\subsection{The finite rule}

We now express the raising coefficients through hook histories, using an
auxiliary number $r$ of rows.

\begin{proposition}\label{prop:finite-rule}
Let $k\geq1$, let $\lambda$ be a partition, and let $r\geq\ell(\lambda)+k$.
Then
\begin{equation}\label{eq:finite-rule}
 b^{(k)}_{\lambda\mu}
 =\sum_{a=0}^{k-1}(-1)^aM^r_{k,a}(\lambda,\mu)-r\beta^k\delta_{\lambda\mu} .
\end{equation}
In particular the right-hand side does not depend on $r$ in this range,
although the individual terms $M^r_{k,a}(\lambda,\mu)$ and
$r\beta^k\delta_{\lambda\mu}$ do.
\end{proposition}

\begin{proof}
By Proposition~\ref{prop:operator-action} with $f=p_k$ and $N=k$,
$b^{(k)}_{\lambda\mu}=[\mu]\Phi_{r+1,r}(p_k)\lambda$. By~\eqref{eq:superpower}
and~\eqref{eq:hook-expansion},
$\Phi_{r+1,r}(p_k)=\sum_a(-1)^as_{H_{k,a}}(\bfU_{r+1})-r\beta^k$. Now apply
the definition of $M^r_{k,a}$. The left-hand side is independent of $r$.
\end{proof}

\subsection{First consequences}
The finite rule gives a preliminary support bound and
homogeneity, followed by the diagonal and top-degree formulas.
The exact off-diagonal support will be determined in
Section~\ref{subsec:compact-rule}.

\begin{corollary}\label{cor:support}
$b^{(k)}_{\lambda\mu}\in\ZZ[\alpha,\beta]$, and $b^{(k)}_{\lambda\mu}=0$ unless
$\lambda\subseteq\mu$ and $|\mu/\lambda|\leq k$. If nonzero, it is homogeneous
of degree $k-|\mu/\lambda|$ in $(\alpha,\beta)$. Hence
\[
 p_kg_\lambda=\sum_{\substack{\mu\supseteq\lambda\\|\mu/\lambda|\leq k}}
 b^{(k)}_{\lambda\mu}\,g_\mu.
\]
\end{corollary}

\begin{proof}
All weights in~\eqref{eq:finite-rule} are integral monomials. By
Lemma~\ref{lem:histories}, a history of a word of length $k$ ending at $\mu$
exists only if $\lambda\subseteq\mu$ and $|\mu/\lambda|\leq k$, and it has
exactly $k-|\mu/\lambda|$ stays. The correction term $r\beta^k\delta_{\lambda\mu}$ has degree $k$ and
occurs only for $\mu=\lambda$.
\end{proof}

The diagonal coefficients have a simple closed form.

\begin{proposition}\label{prop:diagonal}
For every partition $\lambda$ and every $k\geq1$,
\[
 b^{(k)}_{\lambda\lambda}=d(\lambda)\bigl((-\alpha)^k-\beta^k\bigr).
\]
\end{proposition}

\begin{proof}
Let $r\geq\ell(\lambda)+k$. By Lemma~\ref{lem:triangular}, every $X$ in the
algebra generated by $U_1,\dots,U_{r+1}$ satisfies $[\nu]X\kappa=0$ unless
$\kappa\subseteq\nu$. Hence
$[\lambda]XY\lambda=\sum_\nu[\lambda]X\nu\,[\nu]Y\lambda=[\lambda]X\lambda\,[\lambda]Y\lambda$,
so $X\mapsto[\lambda]X\lambda$ is a ring homomorphism sending $U_i$ to
$c_i(\lambda)$. Composed with $\Phi_{r+1}$, it is the evaluation of $\Lam_R$
at $c_1(\lambda),\dots,c_{r+1}(\lambda)$, so
$[\lambda]p_k(\bfU_{r+1})\lambda=\sum_{i=1}^{r+1}c_i(\lambda)^k$. Then exactly $d(\lambda)$ indices $i\leq r+1$ lie in
$\Add^+(\lambda)$, where $c_i(\lambda)=-\alpha$. The index $i=\ell(\lambda)+1$
has $c_i(\lambda)=0$, and the remaining $r-d(\lambda)$ indices have
$\lambda_i=\lambda_{i-1}$ and $c_i(\lambda)=\beta$. By
Proposition~\ref{prop:operator-action} with $f=p_k$
and~\eqref{eq:superpower},
\[
 b^{(k)}_{\lambda\lambda}=[\lambda]p_k(\bfU_{r+1})\lambda-r\beta^k
 =d(\lambda)(-\alpha)^k+(r-d(\lambda))\beta^k-r\beta^k
 =d(\lambda)\bigl((-\alpha)^k-\beta^k\bigr).\qedhere
\]
\end{proof}

For $k=1$ the finite rule reduces to a Pieri rule.

\begin{corollary}[The case $k=1$; {\cite[Prop.~8.12]{Yeliussizov}}]\label{cor:p1}
$p_1g_\lambda=\sum_{\mu}g_\mu-(\alpha+\beta)\,d(\lambda)\,g_\lambda$, the sum over
all $\mu$ obtained from $\lambda$ by adding one box.
\end{corollary}

\begin{proof}
By Proposition~\ref{prop:finite-rule}, for $k=1$ the hook words are the
single letters $i\in[1,r+1]$. The coefficient of
$\lambda+\square_i$ in $U_i\lambda$ is $1$, and no other letter produces
$\lambda+\square_i$. The diagonal is given by
Proposition~\ref{prop:diagonal}.
\end{proof}

The coefficients of top degree are those of the classical rule.

\begin{corollary}\label{cor:top}
If $|\mu/\lambda|=k$, then $b^{(k)}_{\lambda\mu}=(-1)^{\hgt(\mu/\lambda)}$ if
$\mu/\lambda$ is a ribbon, and $b^{(k)}_{\lambda\mu}=0$ otherwise.
\end{corollary}

\begin{proof}
By Corollary~\ref{cor:support}, $b^{(k)}_{\lambda\mu}$ is homogeneous of
degree $0$, hence an integer, when $|\mu/\lambda|=k$. It may therefore be
computed at $\alpha=\beta=0$, where $g^{(0,0)}_\nu=s_\nu$ by~\eqref{eq:duality}
and~\ref{G2}. The claim is then the classical rule~\eqref{eq:classical-MN}.
\end{proof}

Finally, the involution $\omega$ exchanges the roles of $\alpha$ and $\beta$.

\begin{corollary}\label{cor:omega}
$b^{(k)}_{\lambda'\mu'}(\beta,\alpha)=(-1)^{k-1}\,b^{(k)}_{\lambda\mu}(\alpha,\beta)$.
\end{corollary}

\begin{proof}
Apply the ring automorphism $\omega$ to
$p_kg^{(\alpha,\beta)}_\lambda=\sum_\mu b^{(k)}_{\lambda\mu}(\alpha,\beta)g^{(\alpha,\beta)}_\mu$.
By Lemma~\ref{lem:omega-g} and $\omega p_k=(-1)^{k-1}p_k$,
\[
 (-1)^{k-1}p_kg^{(\beta,\alpha)}_{\lambda'}
 =\sum_\mu b^{(k)}_{\lambda\mu}(\alpha,\beta)\,g^{(\beta,\alpha)}_{\mu'} .
\]
The coefficients of $p_kg^{(\beta,\alpha)}_{\lambda'}$ in the basis
$\{g^{(\beta,\alpha)}_\nu\}$ are, by definition, the numbers
$b^{(k)}_{\lambda'\nu}(\beta,\alpha)$; comparing coefficients gives the claim.
\end{proof}

\section{Cancellation and decorated broken ribbons}\label{sec:cancellation}
The finite rule of Proposition~\ref{prop:finite-rule} is alternating and
depends on the auxiliary number $r$. It has two sources of redundancy: disconnected word
supports and addition-free histories. In this section we remove these by sign-reversing
involutions, leaving a finite signed model independent of the cutoff.

\subsection{Statement} We first state the cutoff-free signed rule and then construct the
involutions that remove disconnected supports and addition-free terms.

For a hook word $w$ and a history $h$ of $w$ put
$\swt(w,h)=(-1)^{\leg(w)}\wt(h)$.

\begin{definition}\label{def:Cplus}
For partitions $\lambda,\mu$ and $k\geq1$, let $\cC^+_k(\lambda,\mu)$ be the
set of pairs $(w,h)$ such that
\begin{enumerate}[label=\textup{(\roman*)}]
\item $w$ is a hook word of length $k$ (with arbitrary positive letters)
whose support $\supp(w)$ is an interval of integers;
\item $h$ is a history of $w$ from $\lambda$ with $\enD(h)=\mu$;
\item $h$ contains at least one addition.
\end{enumerate}
\end{definition}

We first check that this set is finite.

\begin{lemma}\label{lem:row-bound}
If $(w,h)\in\cC^+_k(\lambda,\mu)$, then $\max\supp(w)\leq\ell(\lambda)+k$. In
particular $\cC^+_k(\lambda,\mu)$ is finite.
\end{lemma}

\begin{proof}
We write $\supp(w)=[p,q]$ and $\ell=\ell(\lambda)$, and choose a row $j$
in which the history adds a box. If $j>\ell$, the history must add a
box in every row $\ell+1,\ldots,j$. Since $p \leq j \leq q$, this history also uses each of the letters
$j+1,\ldots,q$, so
\[
 k\geq(j-\ell)+(q-j)=q-\ell.
\]
If $j\leq\ell$, the chosen addition together with these latter letters
gives $k\geq1+q-j$, hence $q\leq\ell+k-1$. In either case
$q\leq\ell+k$. There are only finitely many words of length $k$ in
this alphabet, and each has finitely many histories.
\end{proof}

The following theorem is the connected-support form of
Proposition~\ref{prop:finite-rule}.

\begin{theorem}\label{thm:connected-rule}
For all partitions $\lambda,\mu$ and all $k\geq1$,
\[
 b^{(k)}_{\lambda\mu}
 =\sum_{(w,h)\in\cC^+_k(\lambda,\mu)}\swt(w,h)
 +\delta_{\lambda\mu}\,d(\lambda)\bigl((-\alpha)^k-\beta^k\bigr).
\]
\end{theorem}

Theorem~\ref{thm:B} of the introduction is
Theorem~\ref{thm:connected-rule} together with Lemma~\ref{lem:row-bound} and
Proposition~\ref{prop:broken-ribbon} below. For $\lambda=\mu$ the sum over
$\cC^+_k(\lambda,\lambda)$ is empty, since a history with an addition cannot
end at its starting point, and Theorem~\ref{thm:connected-rule} reduces to
Proposition~\ref{prop:diagonal}. The proof of this theorem is given in Subsection~\ref{proof of}.

\subsection{Reduction to interval supports}

The first step removes the hook words whose letters do not form an
interval. This is the cancellation of Fomin and Greene in the proof of
their noncommutative Murnaghan--Nakayama rule \cite[Thm.~1.3]{FominGreene};
the same device is used in \cite{BSZ} and \cite[Lemma~5.3]{NguyenKk}. Since
we need it under the sole hypothesis of far commutativity, we include the
proof.

\begin{lemma}\label{lem:FG}
Let $u_1,\dots,u_n$ be elements of an associative algebra such that
$u_iu_j=u_ju_i$ whenever $|i-j|\geq2$, and put $u_w=u_{w_1}\cdots u_{w_N}$
for a word $w=w_1\cdots w_N$. Then for every $k\geq1$
\[
 \sum_{a=0}^{k-1}(-1)^a\sum_{T\in\SSYT_n(H_{k,a})}u_{w(T)}
 =\sum_w(-1)^{\leg(w)}u_w ,
\]
where the sum on the right runs over the hook words $w$ of length $k$ with
letters in $[n]$ whose support is an interval.
\end{lemma}

\begin{proof}
By Lemma~\ref{lem:hook-words}, the left-hand side is
$\sum_w(-1)^{\leg(w)}u_w$ over all hook words $w$ of length $k$ with letters
in $[n]$. Let $w$ be such a word whose support is not an interval.
Decompose $\supp(w)=I_1\sqcup I_2\sqcup\cdots$ into maximal intervals with
$I_1<I_2<\cdots$, and let $x=\min I_2$. Then $c_1=\min\supp(w)\in I_1$, so
$c_1<x$, and $x-1\notin\supp(w)$. Define $\tau w$ as follows: if
$x\in\Dec(w)$, remove $x$ from $\Dec(w)$ and add one copy of $x$ to
$\Inc(w)$; if $x\notin\Dec(w)$, then $x\in\Inc(w)$ and $x\neq c_1$, and we
remove one copy of $x$ from $\Inc(w)$ and add $x$ to $\Dec(w)$.

In the first case the new decreasing part consists of letters $>c_1$ and the
new increasing part still has minimum $c_1<x$. In the second case the
increasing part keeps the minimum $c_1$, and $\Dec(w)\cup\{x\}$ consists of
letters $>c_1$. So $\tau w$ is a hook word of length $k$ with the same
multiset of letters and the same support as $w$, and
$\leg(\tau w)=\leg(w)\pm1$. Since the support is unchanged, $\tau$ uses the
same letter $x$ for $\tau w$, and it moves $x$ back; as a hook word is
determined by its decreasing part and the multiset of its increasing part,
$\tau(\tau w)=w$.

Finally, if $x\in\Dec(w)$, write
$w=D_{>x}\,x\,D_{<x}\,C_{<x}\,C_{=x}\,C_{>x}$, where $D_{>x}$ and $D_{<x}$ are
the letters of $\Dec(w)$ that are $>x$ and $<x$ in decreasing order, and
$C_{<x},C_{=x},C_{>x}$ are the letters of $\Inc(w)$ that are $<x$, $=x$, $>x$
in increasing order. Then
$\tau w=D_{>x}\,D_{<x}\,C_{<x}\,x\,C_{=x}\,C_{>x}$. All letters of the block
$D_{<x}C_{<x}$ are $<x$, hence lie in $I_1$, hence are $\leq x-2$. So $u_x$
commutes with the corresponding product of generators, and
$u_{\tau w}=u_w$. The case $x\notin\Dec(w)$ is the inverse operation.
Therefore the words whose support is not an interval cancel in pairs.
\end{proof}

By Lemma~\ref{lem:local}, the operators $U_1,\dots,U_n$ satisfy the
hypothesis of Lemma~\ref{lem:FG}.

\subsection{The signed set}

We realize $b^{(k)}_{\lambda\mu}$ as a signed count and then construct a
sign-reversing involution.

Fix $k\geq1$, a partition $\lambda$ and an integer $r\geq\ell(\lambda)+k$. Let
$\cH$ be the set of pairs $(w,h)$ where $w$ is a hook word of length $k$
with letters in $[1,r+1]$ whose support is an interval, and $h$ is a
history of $w$ from $\lambda$. Let $\mathcal V=\{v_2,\dots,v_{r+1}\}$ be a set
of $r$ formal symbols, and put $\enD(v_i)=\lambda$ and $\swt(v_i)=-\beta^k$.
By Proposition~\ref{prop:finite-rule}, Lemma~\ref{lem:hook-words},
Lemma~\ref{lem:FG} with $u_i=U_i$, and Lemma~\ref{lem:histories}, for every
$\mu$
\begin{equation}\label{eq:signed-set}
 b^{(k)}_{\lambda\mu}=\sum_{\substack{x\in\cH\sqcup\mathcal V\\ \enD(x)=\mu}}\swt(x).
\end{equation}
We decompose $\cH=\cH_{\rm free}\sqcup\cH_{\rm sing}\sqcup\cH_{\rm con}$, where
$(w,h)$ lies in $\cH_{\rm free}$ if $|\supp(w)|\geq2$ and $h$ has no
addition, in $\cH_{\rm sing}$ if $|\supp(w)|=1$ and $h$ has no addition, and
in $\cH_{\rm con}$ if $h$ has at least one addition. We construct an
involution $\iota$ of $\cH\sqcup\mathcal V$ that preserves $\enD$, reverses
the sign of $\swt$ on non-fixed points, and whose fixed points we
determine.

\subsection{Addition-free histories}

We next cancel the histories that consist of stays only, for words with at
least two distinct letters.

By Lemma~\ref{lem:triangular}, a word $u$ has a history from $\lambda$
without additions if and only if $c_j(\lambda)\neq0$ for every letter $j$ of
$u$. In that case this history is unique: every state equals $\lambda$, and
its weight is $\prod_jc_{u_j}(\lambda)$. In particular, the existence and the
weight of an addition-free history depend only on the multiset of letters
of $u$.

\begin{lemma}\label{lem:tau-prime}
Let $w$ be a hook word whose support is an interval with at least two
elements, and let $x=\max\supp(w)$. Define $\tau'w$ as follows: if
$x\in\Dec(w)$, then $x=d_a$; delete it and append $x$ at the end of the
increasing part. If $x\notin\Dec(w)$, then $x=c_t$; delete it and put $x$ at
the front of the decreasing part. Then $\tau'w$ is a hook word with the same
multiset of letters and $\leg(\tau'w)=\leg(w)\pm1$, and $\tau'(\tau'w)=w$.
\end{lemma}

\begin{proof}
If $x\in\Dec(w)$, then $x=d_a$ since $x$ is the largest letter. The word
$d_{a-1}\cdots d_1\,c_1\cdots c_tx$ is a hook word, because $x\geq c_t$. If
$x\notin\Dec(w)$, then $x$ occurs in the increasing part, so $c_t=x$. Since
the support has at least two elements, $c_1<x$, so $t\geq2$. Then
$x\,d_a\cdots d_1\,c_1\cdots c_{t-1}$ is a hook word: $x>d_a$ because $x$ is
the maximum and $x\notin\Dec(w)$, and the increasing part $c_1\cdots c_{t-1}$
is nonempty with minimum $c_1$. In the first case $x\notin\Dec(\tau'w)$ and
$x$ is the last letter of $\tau'w$; in the second case $x\in\Dec(\tau'w)$.
So $\tau'$ undoes itself.
\end{proof}

For $(w,h)\in\cH_{\rm free}$ let $\iota(w,h)=(\tau'w,h'')$, where $h''$ is the
unique addition-free history of $\tau'w$. It exists because $\tau'$
preserves the multiset of letters. Then $\iota(w,h)\in\cH_{\rm free}$,
$\iota^2=\mathrm{id}$, $\enD=\lambda$ is preserved, and
$\swt(\iota(w,h))=-\swt(w,h)$.

\subsection{Singletons and the supersymmetric correction}

In this subsection we treat words with a single letter value. Let
$(w,h)\in\cH_{\rm sing}$ and write $\supp(w)=\{i\}$. Every letter of $\Dec(w)$
exceeds $c_1=i$, so $\Dec(w)=\varnothing$, $w=i^k$ and $\leg(w)=0$. The
addition-free history of $i^k$ exists if and only if $c_i(\lambda)\neq0$, and
then $\swt=c_i(\lambda)^k$. For $1\leq i\leq r+1$ we have:
\begin{itemize}
\item $c_i(\lambda)=-\alpha$ if and only if $i\in\Add^+(\lambda)$;
\item $c_i(\lambda)=0$ if and only if $i=\ell(\lambda)+1$;
\item $c_i(\lambda)=\beta$ if and only if $i\in E:=\{i\in[2,r+1]:\lambda_i=\lambda_{i-1}\}$.
\end{itemize}
(For $i=1$ the third case cannot occur since $\lambda_0=\infty$.) Let
$D=[2,r+1]\setminus E=\{i\in[2,r+1]:\lambda_i<\lambda_{i-1}\}$. For $i\in E$,
define $\iota$ to exchange $(i^k,h_0)$, the element with the addition-free
history $h_0$, with $v_i$; their signed weights $\beta^k$ and $-\beta^k$ are
opposite, and both end at $\lambda$. All remaining elements of
$\cH_{\rm sing}\sqcup\mathcal V$ are fixed by $\iota$, namely
\begin{itemize}
\item $(i^k,h_0)$ for $i\in\Add^+(\lambda)$, with signed weight $(-\alpha)^k$;
\item $v_i$ for $i\in D$, with signed weight $-\beta^k$.
\end{itemize}
Since $\lambda_{j+1}<\lambda_j$ forces $\lambda_j>0$ and hence
$j\leq\ell(\lambda)\leq r$, the map $i\mapsto i-1$ is a bijection from $D$
onto $\{j\geq1:\lambda_{j+1}<\lambda_j\}$. Hence $|D|=d(\lambda)$, and
$|\Add^+(\lambda)|=d(\lambda)$. Finally, $\iota$ fixes every element of $\cH_{\rm con}$.

\subsection{Proof of Theorem~\ref{thm:connected-rule}}\label{proof of}

The map $\iota$ is an involution of $\cH\sqcup\mathcal V$ that preserves
$\enD$ and reverses $\swt$ on non-fixed points. Hence~\eqref{eq:signed-set}
reduces to the sum over the fixed points:
\[
 b^{(k)}_{\lambda\mu}
 =\sum_{\substack{(w,h)\in\cH_{\rm con}\\ \enD(h)=\mu}}\swt(w,h)
 +\delta_{\lambda\mu}\bigl(d(\lambda)(-\alpha)^k-d(\lambda)\beta^k\bigr).
\]
The elements of $\cH_{\rm con}$ ending at $\mu$ are exactly the elements of
$\cC^+_k(\lambda,\mu)$ whose letters lie in $[1,r+1]$. By
Lemma~\ref{lem:row-bound} this is all of $\cC^+_k(\lambda,\mu)$, because
$\ell(\lambda)+k\leq r$. \qed

\begin{remark}[Attribution]\label{rem:FG}
Disconnected words are removed by the Fomin--Greene cancellation
(Lemma~\ref{lem:FG}); the same device is used for $k$-Schur functions in
\cite{BSZ} and for $K$-$k$-Schur functions in \cite{NguyenKk}. The
involution on $\cH_{\rm free}$ is needed because the stays of $U_i$ at a
fixed partition generally carry nonzero weights. The pairing of
$\cH_{\rm sing}$ with $\mathcal V$ absorbs the supersymmetric correction;
this pairing is what removes the cutoff $r$.
\end{remark}

\subsection{Decorated histories}

The surviving pairs $(w,h)$ have a simple description, which shows that the
boxes added form a broken ribbon.

\begin{definition}\label{def:decorated}
\emph{Decoration data} of length $k$ consist of integers $1\leq p\leq q$,
integers $f_i\geq0$ and $\epsilon_i\in\{0,1\}$ for $i\in[p,q]$, such that
\[
 f_p\geq1,\qquad \epsilon_p=0,\qquad f_i+\epsilon_i\geq1\ (i\in[p,q]),
 \qquad \sum_{i=p}^q(f_i+\epsilon_i)=k .
\]
The associated \emph{row sequence} is the sequence of rows
\[
 \varrho=\bigl(\underbrace{q,\dots,q}_{f_q},\ \dots,\ \underbrace{p,\dots,p}_{f_p},\
 (p+1)^{\epsilon_{p+1}},\ \dots,\ q^{\epsilon_q}\bigr),
\]
where $i^{\epsilon}$ means that $i$ is present if $\epsilon=1$ and absent if
$\epsilon=0$. A \emph{decorated history} $h$ from $\lambda$ consists of decoration
data together with an admissible move of $U_i$ for each entry $i$ of
$\varrho$, applied in the order of $\varrho$ starting from $\lambda$, with at
least one addition. The first $\sum f_i$ entries form the \emph{first
pass}, and the remaining ones the \emph{second pass}. Its signed weight is $\swt(h)=(-1)^{\sum_i\epsilon_i}\wt(h)$, where $\wt(h)$ is the product of the weights of its moves. Also, we put $\supp(h):=[p,q]$, and we denote by $\DHist_k(\lambda,\mu)$ the set of all decorated histories of length $k$ from $\lambda$ to $\mu$.
\end{definition}

\begin{proposition}\label{prop:broken-ribbon}
The map sending $(w,h)\in\cC^+_k(\lambda,\mu)$ to the data
$\supp(w)=[p,q]$, $f_i=$ multiplicity of $i$ in $\Inc(w)$,
$\epsilon_i=\ind_{\{i\in\Dec(w)\}}$, together with the moves of $h$, is a
bijection onto the set of decorated histories of length $k$ from $\lambda$ to
$\mu$. It preserves signed weights. If $\cC^+_k(\lambda,\mu)\neq\varnothing$,
then $\mu/\lambda$ is a broken ribbon.
\end{proposition}

\begin{proof}
Let $(w,h)\in\cC^+_k(\lambda,\mu)$ with $\supp(w)=[p,q]$. Then $c_1=p$, so
$f_p\geq1$. Every letter of $\Dec(w)$ exceeds $p$, so $\epsilon_p=0$. Every
$i\in[p,q]$ occurs in $w$, so $f_i+\epsilon_i\geq1$, and the total length is
$k$. The word is
\[
 w=q^{\epsilon_q}\cdots(p+1)^{\epsilon_{p+1}}\;p^{f_p}(p+1)^{f_{p+1}}\cdots q^{f_q},
\]
and by~\eqref{eq:two-pass} its letters act in the order of $\varrho$.
Conversely, decoration data determine this word, which is a hook word with
support $[p,q]$ and leg $\sum\epsilon_i$. Histories correspond identically.

Suppose now that some $(w,h)$ exists but $\mu/\lambda$ contains a $2\times2$
block. So there is an $i$ with
$\mu_{i+1}\geq\lambda_i+2$; in particular $i+1\in[p,q]$ because row $i+1$
changes. Row $i+1$ receives at most one action in the second pass, so its
first-pass block raises it from $\lambda_{i+1}$ to at least
$\mu_{i+1}-1\geq\lambda_i+1$. During the first-pass block of row $i+1$, row
$i$ has not yet been acted on: its own first-pass block comes later in
$\varrho$, or does not exist if $i<p$. So row $i$ has length $\lambda_i$ at
that time. At the moment when row $i+1$ would grow from $\lambda_i$ to
$\lambda_i+1$, the rows $i$ and $i+1$ have equal length, and no addition in
row $i+1$ is admissible, a contradiction.
\end{proof}

Geometrically, the additions of a decorated history fill the broken ribbon
$\mu/\lambda$ row by row, from the bottom row $q$ up to row $p$, followed by
at most one further action in each row $i>p$. The stays are \emph{ghosts}
of weight $-\alpha$ or $\beta$, placed on the boundary of the current
partition. Section~\ref{sec:boundary} sums over the ghosts.

\section{The Murnaghan--Nakayama rule}\label{sec:boundary}
We now sum the stay weights row by row. The two-pass order makes the
row contributions independent, giving an interval formula; a further
evaluation reduces it to the single-product rule of
Theorem~\ref{thm:compact-boundary}.

\subsection{Local row factors}

We introduce the generating functions that record first-pass
stays and the possible second-pass move in each row.

Fix partitions $\lambda\subsetneq\mu$ and an integer $k\geq|\mu/\lambda|$,
and put
\[
 n_i=\mu_i-\lambda_i,\qquad e=k-|\mu/\lambda|,\qquad S=\{i\geq1:n_i>0\}.
\]
Let $\cI_e(S)$ be the finite set of intervals $[p,q]$ of positive integers
with $S\subseteq[p,q]$ and $|[p,q]\setminus S|\leq e$. For a row $i\geq1$ and
an integer $y$ define a power series in an auxiliary variable $z$ by
\begin{equation}\label{eq:Phi-def}
 \Phi_i(y)=
 \begin{cases}
  (1+\alpha z)^{-A_i(y)}(1-\beta z)^{-B_i(y)},&\lambda_i\leq y\leq\lambda_{i-1},\\
  0,&\text{otherwise},
 \end{cases}
\end{equation}
where $A_i(y)$ is the number of integers $c$ with
$\max(\lambda_i,1)\leq c\leq y$ and $c<\lambda_{i-1}$, and
$B_i(y)=\ind_{\{i\geq2,\ y=\lambda_{i-1}\}}$. Put
\begin{equation}\label{eq:sigma-def}
 \sigma_i(\mu)=
 \begin{cases}
  \beta,&\mu_i=\mu_{i-1},\\
  -\alpha,&0<\mu_i<\mu_{i-1},\\
  0,&0=\mu_i<\mu_{i-1},
 \end{cases}
\end{equation}
and define
\begin{align}
 P_i(z)&=\Phi_i(\mu_i)-\ind_{\{n_i=0\}},\label{eq:P-def}\\
 R_i(z)&=P_i(z)-z\,\sigma_i(\mu)\,\Phi_i(\mu_i)-\ind_{\{n_i>0\}}\,\Phi_i(\mu_i-1).
 \label{eq:R-def}
\end{align}
These are power series in $z$ with coefficients in $\ZZ[\alpha,\beta]$, and
$P_i$ and $R_i$ depend only on $(\lambda_{i-1},\lambda_i,\mu_{i-1},\mu_i)$. For
$\theta\in\{0,1\}$ we write $\Phi^{[\theta]}_i=\Phi_i(\mu_i-\theta)$, so that
\[
 P_i=\Phi^{[0]}_i-\ind_{\{n_i=0\}},\qquad
 R_i=P_i-z\,\sigma_i(\mu)\,\Phi^{[0]}_i-\ind_{\{n_i>0\}}\,\Phi^{[1]}_i .
\]
Note that $\sigma_i(\mu)=c_i(\mu)$ in the notation of~\eqref{eq:stay-scalar},
and that $\Phi_i(y)$ is nonzero exactly when $\lambda_i\leq y\leq\lambda_{i-1}$;
we then call $y$ \emph{feasible} for row $i$. If $n_i=0$, then $y=\mu_i=\lambda_i$
is feasible, so $\Phi^{[0]}_i(0)=1$, and $P_i$ and $R_i$ have zero constant
term. The set $\cI_e(S)$ is finite because an interval in it has at most
$|S|+e$ elements and contains $\min S$.

\subsection{Timing and the first-pass generating function}
The row above a given action is unchanged during the first pass and
already final during the second. This observation determines the local
generating functions used below.

\begin{lemma}\label{lem:timing}
Consider a decorated history from $\lambda$ with data
$(p,q,(f_i),(\epsilon_i))$ and row sequence $\varrho$.
\begin{enumerate}[label=\textup{(\roman*)}]
\item The rows outside $[p,q]$ are never acted on.
\item For $i\in[p,q]$, the actions on row $i$ are its $f_i$ first-pass
actions, followed by one second-pass action if $\epsilon_i=1$.
\item During the first-pass block of row $i$, row $i-1$ has length
$\lambda_{i-1}$.
\item At the second-pass action of row $i$, row $i-1$ has already reached
its final length.
\end{enumerate}
\end{lemma}

\begin{proof}
(i) and (ii) are read off from $\varrho$. (iii) If $i-1<p$, row $i-1$ is
never acted on. Otherwise its first-pass block comes after that of row $i$,
and all second-pass actions come after all first-pass blocks. (iv) The
first-pass block of row $i-1$ precedes the second pass, and the
second-pass action of row $i-1$, if any, precedes that of row $i$. So all
actions on row $i-1$ precede the second-pass action of row $i$.
\end{proof}

\begin{lemma}\label{lem:first-pass}
Fix a row $i\geq1$ and an integer $y$. Let $\cB_i(y)$ be the set of finite
sequences of moves of $U_i$, applied repeatedly to a partition $\nu$ with
$\nu_{i-1}=\lambda_{i-1}$ fixed and $\nu_i=\lambda_i$ initially, each move
admissible, and ending with $\nu_i=y$. The empty sequence is allowed.
Then
\[
 \sum_{b\in\cB_i(y)}\wt(b)\,z^{\#\mathrm{stays}(b)}=\Phi_i(y).
\]
\end{lemma}

\begin{proof}
Write $L=\lambda_{i-1}$. If $y<\lambda_i$ or $y>L$, then $\cB_i(y)=\varnothing$,
since additions are possible only from levels $c<L$. Let
$\lambda_i\leq y\leq L$. A sequence in $\cB_i(y)$ passes through the levels
$\lambda_i,\lambda_i+1,\dots,y$ and makes exactly one addition from each level
$c<y$. It is therefore determined by the numbers $u_c\geq0$ of stays made at
each level $c\in[\lambda_i,y]$. Every addition in the sequence is
admissible because $c<y\leq L$. By~\eqref{eq:stay-scalar}, a stay at level
$c$ is admissible with weight $-\alpha$ if $1\leq c<L$ and with weight
$\beta$ if $c=L$, and it is not admissible if $c=0<L$. The level $c=L$ can
occur only as $c=y=L$, and then only if $i\geq2$. Summing the geometric
series over the $u_c$ gives a factor $(1+\alpha z)^{-1}$ for each level
counted by $A_i(y)$ and a factor $(1-\beta z)^{-1}$ if $B_i(y)=1$.
\end{proof}

\subsection{The interval formula}
Lemma~\ref{lem:timing} makes the row choices independent once
the initial and final partitions are fixed. Their signed
generating functions therefore multiply for each support
interval.

\begin{theorem}\label{thm:boundary}
For every partition $\lambda$ and every $k\geq1$,
\[
 p_k\,g_\lambda=d(\lambda)\bigl((-\alpha)^k-\beta^k\bigr)\,g_\lambda
 +\sum_{\mu}\Bigl([z^{\,k-|\mu/\lambda|}]\sum_{[p,q]\in\cI_{k-|\mu/\lambda|}(S)}
 P_p(z)\prod_{i=p+1}^{q}R_i(z)\Bigr)\,g_\mu ,
\]
where $\mu$ runs over the partitions $\mu\supsetneq\lambda$ such that
$\mu/\lambda$ is a broken ribbon with $|\mu/\lambda|\leq k$, and $S$,
$P_i$, $R_i$ are formed with respect to $(\lambda,\mu)$. Moreover the
bracketed expression equals $b^{(k)}_{\lambda\mu}$ for every
$\mu\supsetneq\lambda$ with $|\mu/\lambda|\leq k$, and it vanishes if
$\mu/\lambda$ is not a broken ribbon.
\end{theorem}

\begin{proof}
By Corollary~\ref{cor:support},
$p_kg_\lambda=\sum_\mu b^{(k)}_{\lambda\mu}g_\mu$, the sum over
$\mu\supseteq\lambda$ with $|\mu/\lambda|\leq k$, and
$b^{(k)}_{\lambda\lambda}=d(\lambda)((-\alpha)^k-\beta^k)$ by
Proposition~\ref{prop:diagonal}. It therefore suffices to show that for
every $\mu\supsetneq\lambda$ with $e=k-|\mu/\lambda|\geq0$ the coefficient
$b^{(k)}_{\lambda\mu}$ equals the bracketed expression; together with
Proposition~\ref{prop:broken-ribbon}, this also shows that the bracketed
expression vanishes unless $\mu/\lambda$ is a broken ribbon.

Fix such a $\mu$. By Theorem~\ref{thm:connected-rule} and
Proposition~\ref{prop:broken-ribbon}, $b^{(k)}_{\lambda\mu}$ is the signed sum
of all decorated histories of length $k$ from $\lambda$ to $\mu$. We group
them according to $[p,q]$. By
Lemma~\ref{lem:timing}(i) the rows outside $[p,q]$ keep their lengths, so
only intervals with $[p,q]\supseteq S$ occur.

Fix $[p,q]\supseteq S$. To a decorated history with support $[p,q]$ we
associate, for each $i\in[p,q]$, its \emph{row datum}: the sequence $b_i$
of first-pass moves in row $i$, and, if $\epsilon_i=1$, the second-pass
move $m_i$. By Lemma~\ref{lem:timing}(iii), $b_i\in\cB_i(y_i)$, where $y_i$ is
the length of row $i$ after its first pass. By Lemma~\ref{lem:timing}(iv),
the move $m_i$ is applied when rows $i-1$ and $i$ have lengths $\mu_{i-1}$ and
$y_i$. Since row $i$ ends at $\mu_i$, exactly one of the following holds:
\begin{enumerate}[label=(R\arabic*),start=0]
\item $\epsilon_i=0$, and $b_i\in\cB_i(\mu_i)$ is nonempty, since length of $b_i$ is $f_i$ and $f_i+\epsilon_i \geq 1$;
\item $\epsilon_i=1$, $b_i\in\cB_i(\mu_i)$, and $m_i$ is a stay at
$(\mu_{i-1},\mu_i)$ of weight $\sigma_i(\mu)$ (this case occurs only when
$\sigma_i(\mu)\neq0$; if $\sigma_i(\mu)=0$, no admissible stay exists at
$(\mu_{i-1},\mu_i)$ and type~\textup{(R1)} cannot occur);
\item $\epsilon_i=1$, $b_i\in\cB_i(\mu_i-1)$, and $m_i$ is an addition, which
is admissible because $\mu_i-1<\mu_{i-1}$.
\end{enumerate}
Row $p$ is always of type (R0), since $\epsilon_p=0$ and $f_p\geq1$.

Conversely, suppose we choose row data $(b_i,m_i)$ of types
\textup{(R0)}--\textup{(R2)} for all $i\in[p,q]$, with row $p$ of type
\textup{(R0)} and with $e$ stays in total. Set $f_i=|b_i|$, define
$\epsilon_i\in\{0,1\}$ according to the type, and form the row sequence $\varrho=q^{f_q}\cdots p^{f_p}(p+1)^{\epsilon_{p+1}}\cdots q^{\epsilon_q}$.
Then $f_p\geq1$, $\epsilon_p=0$, $f_i+\epsilon_i\geq1$ for all $i\in[p,q]$,
and the total number of moves is $|\mu/\lambda|+e=k$. When we apply the
first-pass sequences $b_q,b_{q-1},\dots,b_p$ in that order starting from
$\lambda$, row $i-1$ has not yet been modified when $b_i$ acts on row $i$, so
row $i-1$ still has length $\lambda_{i-1}$; thus $b_i\in\cB_i(y_i)$ is
admissible and brings row $i$ from $\lambda_i$ to $y_i$. Next, when we apply
the second-pass moves $m_{p+1},\dots,m_q$ from top to bottom, row $i-1$ has
already completed both passes and has length $\mu_{i-1}$ when $m_i$ acts on
row $i$ (at length $y_i$); thus $m_i$ is admissible with the prescribed weight
and brings row $i$ to $\mu_i$. Since $\mu\supsetneq\lambda$, there is at least
one addition, so this defines a decorated history in $\DHist_k(\lambda,\mu)$ of length $k$
with support $[p,q]$, and the two constructions are mutually inverse.

Let
\[
 \DHist(\lambda,\mu)
 :=\bigsqcup_{L\geq|\mu/\lambda|}\DHist_L(\lambda,\mu).
\]
A history in $\DHist_L(\lambda,\mu)$ has exactly
$L-|\mu/\lambda|$ stays. Hence extracting the coefficient of
$z^{k-|\mu/\lambda|}$ from the generating function over all
lengths selects precisely the histories of length $k$:
\[
 b^{(k)}_{\lambda\mu}
 =[z^{k-|\mu/\lambda|}]
 \sum_{h\in\DHist(\lambda,\mu)}
       \swt(h)z^{\#\mathrm{stays}(h)}.
\]
For a fixed number of stays, the length is fixed and the row
bound in Lemma~\ref{lem:row-bound} makes this coefficient a
finite sum.

A decorated history with support $[p,q]$ is determined by its row data
$(b_i,m_i)_{i\in[p,q]}$, where $m_i$ is absent if $\epsilon_i=0$ (so $m_p$ is
always absent). For a row datum $x=(b_i,m_i)$ of row $i$ put
\[
 S_i(x)=(-1)^{\epsilon_i}\,\wt(b_i)\,\wt(m_i)\,
 z^{\#\mathrm{stays}(b_i)+\ind_{\{m_i\text{ is a stay}\}}},
\]
where $\wt(m_i)=1$ if $m_i$ is absent or an addition. Let $T^{(i)}_0,T^{(i)}_1,T^{(i)}_2$
be the sets of row data of row $i$ of types \textup{(R0)}, \textup{(R1)},
\textup{(R2)}, and $T^{(i)}=T^{(i)}_0\sqcup T^{(i)}_1\sqcup T^{(i)}_2$. Each row of $[p,q]\setminus S$ contributes at least one stay.
Thus, modulo $z^{e+1}$, only intervals in $\cI_e(S)$ can
contribute. The following sum and its rowwise factorization
are therefore well defined coefficientwise. Since row $p$
is of type \textup{(R0)} and the rows are independent,
\begin{equation}\label{eq:A}
 b^{(k)}_{\lambda\mu}=[z^{\,k-|\mu/\lambda|}]\sum_{[p,q]\supseteq S}
 \Bigl(\sum_{x\in T^{(p)}_0}S_p(x)\Bigr)\prod_{i=p+1}^{q}\Bigl(\sum_{x\in T^{(i)}}S_i(x)\Bigr).
\end{equation}
By Lemma~\ref{lem:first-pass},
\begin{align}
 \sum_{x\in T^{(i)}_0}S_i(x)&=\Phi^{[0]}_i-\ind_{\{n_i=0\}}=P_i, \label{eq:B}\\
 \sum_{x\in T^{(i)}_1}S_i(x)&=-z\sigma_i(\mu)\Phi^{[0]}_i, \label{eq:C}\\
 \sum_{x\in T^{(i)}_2}S_i(x)&=-\Phi^{[1]}_i. \label{eq:D}
\end{align}
(In \eqref{eq:B} the empty sequence lies in $\cB_i(\mu_i)$ exactly when $\mu_i=\lambda_i$,
with weight $1$ and no stays; \eqref{eq:D} vanishes when $n_i=0$.) Hence
$\sum_{T^{(p)}_0}S_p=P_p$ and, for $p<i\leq q$, $\sum_{T^{(i)}}S_i=P_i-z\sigma_i(\mu)\Phi^{[0]}_i-\Phi^{[1]}_i=R_i$,
so
\[
 b^{(k)}_{\lambda\mu}=[z^{\,k-|\mu/\lambda|}]\sum_{[p,q]\supseteq S}P_p\prod_{i=p+1}^qR_i=[z^{\,k-|\mu/\lambda|}]\sum_{[p,q]\in \cI_e(S) }P_p\prod_{i=p+1}^qR_i,
\]
since only intervals in $\cI_e(S)$ can
contribute.
\end{proof}

By Hall duality, Theorem~\ref{thm:boundary} is also a lowering rule for the
canonical basis.

\begin{corollary}\label{cor:G-lowering}
For every partition $\mu$ and every $k\geq1$,
\[
 p_k^\perp\,G_\mu=\sum_{\lambda\subseteq\mu}b^{(k)}_{\lambda\mu}\,G_\lambda ,
\]
a finite sum, where $p_k^\perp=k\,\partial/\partial p_k$ acts on $\hLam_R$
degreewise.
\end{corollary}

\begin{proof}
By~\ref{G1}, every $F\in\hLam_R$ has a unique expansion
$F=\sum_\lambda c_\lambda G_\lambda$, and $c_\lambda=\Hall{F}{g_\lambda}$: since
$g_\lambda$ has degree $\leq|\lambda|$ and $G_\nu$ has no component of degree
$<|\nu|$, only finitely many terms contribute to the pairing. Take
$F=p_k^\perp G_\mu$. By Lemma~\ref{lem:hopf}(c),
$\Hall{p_k^\perp G_\mu}{g_\lambda}=\Hall{G_\mu}{p_kg_\lambda}=b^{(k)}_{\lambda\mu}$,
which vanishes unless $\lambda\subseteq\mu$ and $|\mu/\lambda|\leq k$ by
Corollary~\ref{cor:support}.
\end{proof}

\subsection{A single-product formula}
\label{subsec:compact-rule}

The interval sum in Theorem~\ref{thm:boundary} admits a further
cancellation. It reduces to one rational generating function and shows
that, apart from the skew diagram, only two endpoint contacts are needed.

\begin{definition}\label{def:corner-connected}
Let $\lambda\subsetneq\mu$. We call $\mu/\lambda$ a
\emph{corner-connected broken ribbon} if its occupied rows form an
interval $[a,b]$ and
\[
 \mu_i\in\{\lambda_{i-1},\lambda_{i-1}+1\}
 \qquad(a<i\leq b).
\]
Equivalently, it contains no $2\times2$ block and is connected when
boxes sharing either an edge or a corner are considered adjacent. Let $D=|\mu/\lambda|$, let $r=b-a+1$, and let $c$ be its number of
edge-connected components. Define
\[
 \varepsilon_{\mathrm{top}}
 =\ind_{\{a>1,\ \mu_a=\lambda_{a-1}\}},\qquad
 \varepsilon_{\mathrm{bot}}
 =\ind_{\{\lambda_b=\lambda_{b+1}\}},
\]
and put
\[
 B=r-1+\varepsilon_{\mathrm{top}}+\varepsilon_{\mathrm{bot}},
 \qquad A=D+c-B.
\]
Here $c$ denotes the number of components, not the number of columns. Both $A$ and $B$ are nonnegative integers. Indeed,
$0\leq B\leq r+1$, while $D\geq r$ and $c\geq1$, so
$D+c\geq r+1$ and hence $A=D+c-B\geq0$.
\end{definition}

\begin{example}\label{ex:corner-connected}
Let $\lambda=(2,1)$ and $\mu=(3,2,2)$. Then $\mu/\lambda$ occupies rows
$[a,b]=[1,3]$ with boxes $(1,3)$, $(2,2)$, $(3,1)$, $(3,2)$. Since
$\mu_2=2=\lambda_1$ and $\mu_3=2=\lambda_2+1$, the skew shape $\mu/\lambda$ is
a corner-connected broken ribbon of size $D=4$ and $r=b-a+1=3$, consisting of
$c=2$ edge-connected components, namely $\{(1,3)\}$ and
$\{(2,2),(3,1),(3,2)\}$ which meet only at the corner between $(1,3)$ and
$(2,2)$. Here $\varepsilon_{\mathrm{top}}=0$ (as $a=1$) and
$\varepsilon_{\mathrm{bot}}=1$ (as $\lambda_3=\lambda_4=0$), giving $B=3$ and
$A=3$.
\end{example}

\begin{theorem}\label{thm:compact-boundary}
Let $\lambda\subsetneq\mu$. If $\mu/\lambda$ is not a corner-connected broken ribbon, then $b^{(k)}_{\lambda\mu}=0$ for
every $k\geq1$. Otherwise, with the statistics of
Definition~\ref{def:corner-connected},
\begin{equation}\label{eq:compact-boundary}
 b^{(k)}_{\lambda\mu}(\alpha,\beta)
 =(-1)^{r-c}(\alpha+\beta)^{c-1}
 [z^{\,k-D-c+1}]
 \frac{1}{(1+\alpha z)^A(1-\beta z)^B}.
\end{equation}
Coefficients with negative extraction index are understood to be zero.
Equivalently,
\begin{equation}\label{eq:compact-generating}
 \sum_{e\geq0}b^{(D+e)}_{\lambda\mu}z^e
 =\frac{(-1)^{r-c}(\alpha+\beta)^{c-1}z^{c-1}}
 {(1+\alpha z)^A(1-\beta z)^B}.
\end{equation}
\end{theorem}

\begin{proof}
We evaluate the factors in Theorem~\ref{thm:boundary}. Write
$S=\{i:\mu_i>\lambda_i\}$, $a=\min S$, and $b=\max S$.
For an active row $i\in S$, formula~\eqref{eq:R-def} becomes
\begin{equation}\label{eq:active-R}
 R_i=(1-z\sigma_i(\mu))\Phi_i(\mu_i)-\Phi_i(\mu_i-1).
\end{equation}
If row $i-1$ is unchanged, then $\mu_{i-1}=\lambda_{i-1}$ and 
$R_i=0$. Indeed, for $\mu_i<\lambda_{i-1}$ (in this case $\sigma_i(\mu)=-\alpha , B_i(\mu_i)=B_i(\mu_i-1)=0$),
$$R_i=(1+\alpha z)(1+ \alpha z)^{-A_i(\mu_i)}- (1+ \alpha z)^{-A_i(\mu_i-1)}.$$
Since $A_i(\mu_i)=A_i(\mu_i-1) +1 $, we have $R_i=0$. The same argument with
$\lambda_0=\mu_0=\infty$ applies when $i=1$.

On the other hand, for
$\mu_i=\lambda_{i-1}$ (in this case $\sigma_i(\mu)=\beta$),
$$ R_i= (1-\beta z)(1+\alpha z)^{-A_i(\mu_i)} (1 - \beta z)^{-1}- (1+\alpha z)^{-A_i(\mu_i-1)}.$$
As $A_i(\mu_i)=A_i(\mu_i-1) $, we obtain $R_i=0$.

For an inactive row, $\mu_i=\lambda_i$, and
\begin{equation}\label{eq:inactive-R}
 R_i=(1-z\sigma_i(\mu))\Phi_i(\lambda_i)-1.
\end{equation}
If row $i-1$ is also unchanged, then
$\sigma_i(\mu)=c_i(\lambda)$ and
\[
 \Phi_i(\lambda_i)=\frac{1}{1-zc_i(\lambda)}.
\]
Indeed, the only first-pass moves returning row $i$ to its
initial length are stays. This includes the case
$c_i(\lambda)=0$, when only the empty sequence is allowed.
Substituting into \eqref{eq:inactive-R} gives $R_i=0$. 
It follows that every interval $[p,q]$ with $p<a$ contributes zero:
its product includes $R_a=0$. Every interval with $q>b+1$ also
contributes zero, since it includes $R_{b+2}=0$. If $S$ has a gap,
the first active row $i$ after a gap gives $R_i=0$ and since $S \subseteq [p,q]$, $[p,q]$ will contribute $0$. Thus all coefficients vanish unless $S=[a,b]$.

Assume henceforth that $S=[a,b]$. Then among all intervals $[p,q]$ containing
$[a,b]$, only $[a,b]$ and $[a,b+1]$ can be nonzero, contributing
$P_a\prod_{i=a+1}^b R_i$ and $\bigl(P_a\prod_{i=a+1}^b R_i\bigr)R_{b+1}$,
respectively. Summing the two contributions in Theorem~\ref{thm:boundary} over
$e\geq0$ gives
\begin{equation}\label{eq:collapsed-intervals}
 \sum_{e\geq0}b^{(D+e)}_{\lambda\mu}z^e
 =P_a\prod_{i=a+1}^b R_i\,(1+R_{b+1}).
\end{equation}
(Note that $R_{b+1}$ has zero constant term, so the $[a,b+1]$ term does not
contribute at $e=0$.)

By definition, $a=\min S$ is the first active row and hence $a-1$ is inactive. Thus we obtain $\lambda_a < \mu_a \leq \mu_{a-1}=\lambda_{a-1} $. Consequently we have
\begin{equation}
    \label{eq:New-0}
P_a=\Phi_a(\mu_a)=
  (1+\alpha z)^{-A_a(\mu_a)}(1-\beta z)^{-B_a(\mu_a)}.
\end{equation}
For $a<i\leq b$, put $L=\lambda_{i-1}$. The preceding row is active,
so $\mu_{i-1}>L$. Directly from \eqref{eq:active-R},
\begin{equation}\label{eq:active-R-cases}
 R_i=
 \begin{cases}
  -\Phi_i(L),&\mu_i=L+1,\\
  (\alpha+\beta)z\,\Phi_i(L),&\mu_i=L,\\
  0,&\text{otherwise}.
 \end{cases}
\end{equation}
For the middle case, $0\leq\lambda_i<L<\mu_{i-1}$, so
$\Phi_i(L-1)=(1-\beta z)\Phi_i(L)$ and $\sigma_i(\mu)=-\alpha$, giving
$(1+\alpha z)\Phi_i(L)-(1-\beta z)\Phi_i(L)=(\alpha+\beta)z\,\Phi_i(L)$. For $\mu_i>L+1$, both
$\Phi_i(\mu_i)$ and $\Phi_i(\mu_i-1)$ are zero. 
For $\mu_i<L$, $0 < \mu_i < \mu_{i-1}$ so $\sigma_i(\mu)=-\alpha$. As $B_i(\mu_i),B_i(\mu_i-1)$ are both equals to $0$, we have
$$ R_i=(1+\alpha z) (1+\alpha z)^{-A_i(\mu_i)}-(1+\alpha z)^{-A_i(\mu_i-1)}. $$
Since $ A_i(\mu_i-1)=A_i(\mu_i)-1$, we conclude that $R_i=0$. This proves the
support condition in the statement.

For the inactive row $b+1$,
\eqref{eq:inactive-R} gives
\begin{align*}
 1+R_{b+1}&=(1-z\sigma_{b+1}(\mu))\Phi_{b+1}(\lambda_{b+1})\\
 &=(1-z\sigma_{b+1}(\mu))(1+\alpha z)^{-A_{b+1}(\lambda_{b+1})}(1-\beta z)^{-B_{b+1}(\lambda_{b+1})}.
\end{align*}
Since $\mu_{b+1}=\lambda_{b+1}\leq\lambda_b<\mu_b$, we have
$A_{b+1}(\lambda_{b+1})=\ind_{\{0<\lambda_{b+1}<\lambda_b\}}$,
$B_{b+1}(\lambda_{b+1})=\ind_{\{\lambda_{b+1}=\lambda_b\}}$, and
$\sigma_{b+1}(\mu)=-\alpha$ if $\lambda_{b+1}>0$ while $\sigma_{b+1}(\mu)=0$ if
$\lambda_{b+1}=0$. Evaluating the three possible cases gives
\begin{equation}
    \label{eq:New}
     1+R_{b+1}=
 \begin{cases}
 (1+\alpha z)/(1-\beta z),&\lambda_b=\lambda_{b+1}>0,\\
 (1-\beta z)^{-1},&\lambda_b=\lambda_{b+1}=0,\\
 1,&\lambda_b>\lambda_{b+1}.
 \end{cases}
\end{equation}
Using \eqref{eq:New-0},\eqref{eq:active-R-cases} and\eqref{eq:New}, we can say that the expression in the right hand side of \eqref{eq:collapsed-intervals} will be non-zero, if $\mu_i \in \{\lambda_{i-1},\lambda_{i-1}+1\}$ for $a < i \leq b$. Also, we have already seen $S=[a,b]$. Therefore we assume that $\mu/\lambda$ is corner-connected. Put
\[
 h=\#\{i\in(a,b]:\mu_i=\lambda_{i-1}+1\},\qquad
 t=\#\{i\in(a,b]:\mu_i=\lambda_{i-1}\}.
\]
An edge connection joins consecutive rows in the first case, whereas
the second case separates two edge-connected components. Consequently
$h=r-c$ and $t=c-1$.

Now we show the exponents of $(1-\beta z)^{-1}$ and $ (1+ \alpha z)^{-1}$ in \eqref{eq:collapsed-intervals} are $B, A$ respectively. The top factor, i.e., $P_a=\Phi_a(\mu_a)$, supplies $\varepsilon_{\mathrm{top}}$ copies of
$(1-\beta z)^{-1}$, each of $R_i$ for $a<i\leq b$ supplies one $(1-\beta z)^{-1}$,
and the last factor, i.e., $ 1+R_{b+1}$, supplies $\varepsilon_{\mathrm{bot}}$ copies.
The total exponent of $(1-\beta z)^{-1}$ in \eqref{eq:collapsed-intervals} is therefore $B$.

For completeness, we also compute the exponent of $(1+\alpha z)^{-1}$.
Write $n_i=\mu_i-\lambda_i$ and, for $a<i\leq b$, put
$\eta_i=\mu_i-\lambda_{i-1}\in\{0,1\}$. We recall that for $\lambda_i \leq y \leq \lambda_{i-1}$, $A_{i}(y)$ is the number of integers $c$ with $c < \lambda_{i-1}$ and $\max(\lambda_i,1) \leq c \leq y$. Then
\begin{align*}
A_a(\mu_a)&=\mu_a-\max(\lambda_a,1)+1-\varepsilon_{\mathrm{top}}
  =n_a+1-\ind_{\{\lambda_a=0\}}-\varepsilon_{\mathrm{top}},\\
 \text{and for } a<i\leq b,\quad
 A_i(\lambda_{i-1})&=\lambda_{i-1}-\lambda_i-\ind_{\{\lambda_i=0<\lambda_{i-1}\}}
  =n_i-\eta_i-\ind_{\{\lambda_i=0<\lambda_{i-1}\}}.
  \end{align*} (The integers counted by $A_i(\lambda_{i-1})$ are those in $[\max(\lambda_i,1),\lambda_{i-1}-1]$: there are
$\lambda_{i-1}-\lambda_i$ of them if $\lambda_i>0$, $\lambda_{i-1}-1$ if $\lambda_i=0<\lambda_{i-1}$, and none if $\lambda_i=\lambda_{i-1}=0$.)

Since $\lambda_a\geq\lambda_{a+1}\geq\cdots\geq\lambda_b\geq0$ is weakly
decreasing, if $\lambda_b=0$ then either $\lambda_a=0$ or $\lambda_i$ drops
from positive to zero at a unique index $i\in(a,b]$; hence
\[
 \ind_{\{\lambda_a=0\}}+
 \sum_{i=a+1}^b\ind_{\{\lambda_i=0<\lambda_{i-1}\}}
 =\ind_{\{\lambda_b=0\}}.
\]
Since $\sum_{i=a+1}^b\eta_i=h$, the active rows give exponent
$D-h+1-\varepsilon_{\mathrm{top}}-\ind_{\{\lambda_b=0\}}$.
The last factor subtracts one exactly when
$\lambda_b=\lambda_{b+1}>0$. Since $\ind_{\{\lambda_b=0\}}+\ind_{\{\lambda_b=\lambda_{b+1}>0\}} =\varepsilon_{\mathrm{bot}}$ and $\lambda_b=0$ forces $\lambda_{b+1}=0$, the resulting exponent is
\[
 D-h+1-\varepsilon_{\mathrm{top}} 
 - \ind_{\{\lambda_b=0\}}-\ind_{\{\lambda_b=\lambda_{b+1}>0\}}
=D-h+1-\varepsilon_{\mathrm{top}} 
 -\varepsilon_{\mathrm{bot}}
 =D+c-B=A.
\]
Substitution in \eqref{eq:collapsed-intervals} proves
\eqref{eq:compact-generating}; coefficient extraction proves
\eqref{eq:compact-boundary}.
\end{proof}
The product formula determines exactly which coefficient
polynomials are nonzero.

\begin{corollary}
\label{cor:sharp-support}
For $\lambda\subsetneq\mu$, the polynomial
$b^{(k)}_{\lambda\mu}(\alpha,\beta)$ is nonzero if and only if
$\mu/\lambda$ is a corner-connected broken ribbon and
\[
 k\geq |\mu/\lambda|+c-1,
\]
where $c$ is its number of edge-connected components. In this case it
is divisible by $(\alpha+\beta)^{c-1}$, and
\[
 b^{(|\mu/\lambda|+c-1)}_{\lambda\mu}
 =(-1)^{r-c}(\alpha+\beta)^{c-1}.
\]
\end{corollary}

\begin{proof}
By Theorem~\ref{thm:compact-boundary}, $b^{(k)}_{\lambda\mu}(\alpha,\beta)$
vanishes unless $\mu/\lambda$ is a corner-connected broken ribbon and
$N:=k-D-c+1\geq0$, and when these conditions hold it is divisible by
$(\alpha+\beta)^{c-1}$. Conversely, for every $N\geq0$,
\[
 [z^N](1+\alpha z)^{-A}(1-\beta z)^{-B}
\]
is a nonzero polynomial: if $A>0$, its $\alpha^N$ coefficient is
$(-1)^N\binom{A+N-1}{N}\neq0$; if $A=0$, then $B=D+c>0$ and its
$\beta^N$ coefficient is $\binom{B+N-1}{N}\neq0$. Since
$\ZZ[\alpha,\beta]$ is an integral domain, $b^{(k)}_{\lambda\mu}(\alpha,\beta)$
is nonzero. Finally, for $N=0$ the constant term of the rational function is
$1$, giving the last formula. (Note that nonvanishing holds in the polynomial
ring $\ZZ[\alpha,\beta]$; particular numerical specializations of
$\alpha,\beta$ may vanish.)
\end{proof}

The factor contributed by the edge-connected components has
the following exact multiplicity.

\begin{corollary}
\label{cor:exact-component-divisibility}
Suppose that $\mu/\lambda$ is a corner-connected broken ribbon with
$D$ boxes, $r$ occupied rows, and $c$ edge-connected components.
For every $k\geq D+c-1$, the polynomial
$b^{(k)}_{\lambda\mu}(\alpha,\beta)$ is divisible by
$(\alpha+\beta)^{c-1}$ but not by $(\alpha+\beta)^c$.
\end{corollary}

\begin{proof}
Put $N=k-D-c+1$. By Theorem~\ref{thm:compact-boundary},
the quotient
\[
 Q_{\lambda\mu}^{(k)}(\alpha,\beta)
 :=\frac{b^{(k)}_{\lambda\mu}(\alpha,\beta)}
          {(\alpha+\beta)^{c-1}}
\]
is a polynomial. Since $A+B=D+c$, we obtain
\begin{align*}
 Q_{\lambda\mu}^{(k)}(\alpha,-\alpha)
 &=(-1)^{r-c}[z^N](1+\alpha z)^{-(D+c)}\\
 &=(-1)^{r-c}(-\alpha)^N
   \binom{D+c+N-1}{N}\\
 &=(-1)^{r-c}(-\alpha)^N\binom{k}{N}.
\end{align*}
This is a nonzero polynomial in $\alpha$. Hence
$Q_{\lambda\mu}^{(k)}$ is not divisible by $\alpha+\beta$,
which proves the assertion.
\end{proof}

For $M,N\geq0$, let
\[
 H_M(N)=[z^N](1-z)^{-M}
 =\begin{cases}
 \binom{M+N-1}{N},&M>0,\\
 \delta_{N0},&M=0.
 \end{cases}
\]
Set $H_M(N)=0$ for $N<0$. Thus $H_M(N)$ counts weak compositions of
$N$ into $M$ parts, including the usual empty-composition convention.

\begin{corollary}\label{cor:compact-binomial}
With the notation of Theorem~\ref{thm:compact-boundary}, put
$N=k-D-c+1$. If $N\geq0$, then
\[
 b^{(k)}_{\lambda\mu}
 =(-1)^{r-c}(\alpha+\beta)^{c-1}
 \sum_{u+v=N}(-\alpha)^u\beta^v H_A(u)H_B(v).
\]
In the usual dual stable normalization, $(\alpha,\beta)=(0,1)$,
\begin{equation}\label{eq:dual-stable-binomial}
 b^{(k)}_{\lambda\mu}(0,1)=(-1)^{r-c}H_B(k-D-c+1).
\end{equation}
The coefficient is zero outside the corner-connected support, and the
diagonal coefficient is $b^{(k)}_{\lambda\lambda}(0,1)=-d(\lambda)$.
\end{corollary}

\begin{proof}
We expand each factor of \eqref{eq:compact-boundary} by the binomial
theorem and then specialize. The diagonal formula is
Proposition~\ref{prop:diagonal}.
\end{proof}

\begin{remark}
For a corner-connected broken ribbon, the statistic $B$ also
has the form
\[
 B=\sum_{i\geq1}
 \left(
   \ind_{\{\mu_i>\lambda_{i+1}\}}
   -\ind_{\{\lambda_i>\mu_{i+1}\}}
 \right).
\]
Indeed, the internal indices $a\leq i<b$ each contribute $1$,
the index $a-1$, when present, contributes
$\varepsilon_{\mathrm{top}}$, and the index $b$ contributes
$\varepsilon_{\mathrm{bot}}$; all other summands vanish.
Consequently, \eqref{eq:dual-stable-binomial} agrees with
\cite[Thm.~2.1]{HiepThuy}.
\end{remark}

\subsection{The coefficients are not functions of the skew diagram}

In the classical rule, and in the rules of \cite{NHST,Kundu} on the
$G$-side, the coefficient of a basis element depends only on the skew
diagram. The following example shows that this fails for the raising
rule: the coefficient depends also on the surrounding parts of $\lambda$
through the boundary indicators $\varepsilon_{\mathrm{top}}$ and
$\varepsilon_{\mathrm{bot}}$.

\begin{example}\label{prop:not-shape-only}
Consider the pairs $(\lambda,\mu)=((2,1),(2,2))$ and
$(\lambda,\mu)=((3,1),(3,2))$. In both cases $\mu/\lambda$ is the single box
$(2,2)$ with $a=b=2$, so in the notation of
Definition~\ref{def:corner-connected} both pairs have $D=r=c=1$ and
$\varepsilon_{\mathrm{bot}}=\ind_{\{\lambda_2=\lambda_3\}}=0$. However, for
$((2,1),(2,2))$ we have $\mu_2=2=\lambda_1$, giving
$\varepsilon_{\mathrm{top}}=1$ and $(A,B)=(1,1)$, whereas for
$((3,1),(3,2))$ we have $\mu_2=2<\lambda_1=3$, giving
$\varepsilon_{\mathrm{top}}=0$ and $(A,B)=(2,0)$. Applying
Theorem~\ref{thm:compact-boundary} with $k=2$ yields
\[
 b^{(2)}_{(2,1),(2,2)}
 =[z]\,\frac{1}{(1+\alpha z)(1-\beta z)}=\beta-\alpha,
 \qquad
 b^{(2)}_{(3,1),(3,2)}
 =[z]\,\frac{1}{(1+\alpha z)^2}=-2\alpha.
\]
In particular, $b^{(k)}_{\lambda\mu}$ is not a function of the skew diagram
$\mu/\lambda$ alone.
\end{example}

\begin{example}\label{ex:disconnected}
Let $\lambda=(2,1)$ and $\mu=(3,2)$. Then $\mu/\lambda$ consists of the two
boxes $(1,3)$ and $(2,2)$, which occupy rows $[a,b]=[1,2]$ and meet only at
a corner (since $\mu_2=2=\lambda_1$). Thus $\mu/\lambda$ is a
corner-connected broken ribbon with $D=2$, $r=2$, and $c=2$ edge-connected
components. Here $\varepsilon_{\mathrm{top}}=0$ (as $a=1$) and
$\varepsilon_{\mathrm{bot}}=0$ (as $\lambda_2=1>\lambda_3=0$), so $B=1$ and
$A=D+c-B=3$. By Theorem~\ref{thm:compact-boundary}, for every $k\geq1$,
\[
 b^{(k)}_{(2,1),(3,2)}
 =(\alpha+\beta)\,[z^{k-3}]\,\frac{1}{(1+\alpha z)^3(1-\beta z)}.
\]
In particular, $b^{(2)}_{(2,1),(3,2)}=0$, while at the threshold degree
$k=D+c-1=3$ (Corollary~\ref{cor:sharp-support}) we obtain
$b^{(3)}_{(2,1),(3,2)}=\alpha+\beta$.
\end{example}

\section{Consequences of the rule}\label{sec:consequences}
In this section, we first derive the hook expansion of power sums and examine
the one-parameter specialization.
We then record the stable $K$-homological interpretation and
the action of deformed power sums.

\subsection{The expansion of \texorpdfstring{$p_k$}{pk}}\label{sec:hooks}

Since $g_\one=1$, taking $\lambda=\one$ expands the power sums themselves.
Starting from the empty partition forces the corner-connected broken ribbon
to be a hook, so Theorem~\ref{thm:compact-boundary} gives the complete
expansion without a separate analysis of histories.

\begin{theorem}\label{thm:empty-partition}
For every $k\geq1$,
\begin{equation}\label{eq:pk-expansion}
 p_k=\sum_{\substack{m\geq1,\ j\geq0\\ m+j\leq k}}c^{(k)}_{m,j}\,g_{(m,1^j)},
 \qquad
 c^{(k)}_{m,j}=(-1)^j\,[t^{\,k-m-j}]\,\frac{1}{(1+\alpha t)^m(1-\beta t)^{j+1}} .
\end{equation}
Explicitly,
$c^{(k)}_{m,j}=(-1)^j\sum_{u+v=k-m-j}(-\alpha)^u\binom{m+u-1}{u}\beta^v\binom{j+v}{v}$.
\end{theorem}

\begin{proof}
Take $\lambda=\one$ in Theorem~\ref{thm:compact-boundary}. A nonempty
straight diagram is a corner-connected broken ribbon exactly when
it is a hook $\mu=(m,1^j)$. For this hook,
\[
 D=m+j,\qquad r=j+1,\qquad c=1,\qquad
 \varepsilon_{\mathrm{top}}=0,\qquad
 \varepsilon_{\mathrm{bot}}=1,
\]
so $A=m$ and $B=j+1$. Formula~\eqref{eq:compact-boundary} gives the
claimed coefficient, and the empty-diagram coefficient is zero by
$d(\one)=0$.
\end{proof}

\begin{example}\label{ex:pk-small}
For $k=2,3$, Theorem~\ref{thm:empty-partition} gives
\begin{align*}
 p_2&=g_{(2)}-g_{(1,1)}+(\beta-\alpha)\,g_{(1)},\\
 p_3&=g_{(3)}-g_{(2,1)}+g_{(1,1,1)}+(\beta-2\alpha)\,g_{(2)}+(\alpha-2\beta)\,g_{(1,1)}
 +(\alpha^2-\alpha\beta+\beta^2)\,g_{(1)} .
\end{align*}
These agree with the Schur expansions $g_{(2)}=s_2+\alpha s_1$ and
$g_{(1,1)}=s_{11}+\beta s_1$ \cite[Prop.~8.9]{Yeliussizov}; for instance
$g_{(2)}-g_{(1,1)}+(\beta-\alpha)g_{(1)}
=s_2+\alpha s_1-s_{11}-\beta s_1+(\beta-\alpha)s_1=s_2-s_{11}=p_2$.
\end{example}

\subsection{One-parameter specialization and reduction}
\label{sec:one-parameter}
At $\alpha=0$ the dual canonical functions become the one-parameter dual
stable Grothendieck functions, and the stays of weight $-\alpha$ disappear. 

Recall $\sg^{(\gamma)}_\lambda=g^{(0,-\gamma)}_\lambda$ from~\eqref{eq:one-parameter-translation}.
At $(\alpha,\beta)=(0,-\gamma)$ the stays of weight $-\alpha$ vanish, so every
contributing history is \emph{deterministic}: the operator $U_i$ adds a box
in row $i$ whenever $\nu_i<\nu_{i-1}$, and otherwise stays with weight
$-\gamma$. For $r\geq0$ let $N^r_{k,a}(\lambda,\mu)$ be the number of
$T\in\SSYT_{r+1}(H_{k,a})$ whose deterministic history from $\lambda$ ends at
$\mu$.

\begin{corollary}\label{cor:one-param-finite}
For $r\geq\ell(\lambda)+k$, \[
 p_k\,\sg^{(\gamma)}_\lambda
 =\sum_{\substack{\mu\supseteq\lambda\\|\mu/\lambda|\leq k}}
 \left(
 (-\gamma)^{k-|\mu/\lambda|}
 \sum_{a=0}^{k-1}(-1)^aN^r_{k,a}(\lambda,\mu)
 -r(-\gamma)^k\delta_{\lambda\mu}
 \right)\sg^{(\gamma)}_\mu.
\]
Equivalently, by Theorem~\ref{thm:connected-rule}, the coefficient of
$\sg^{(\gamma)}_\mu$ for $\mu\neq\lambda$ is $(-\gamma)^{k-|\mu/\lambda|}$ times the
signed number $\sum(-1)^{\leg(w)}$ of hook words $w$ of length $k$ with
interval support whose deterministic history from $\lambda$ ends at $\mu$. The
diagonal coefficient is $-d(\lambda)(-\gamma)^k$.
\end{corollary}

\begin{proof}
The substitution $(\alpha,\beta)\mapsto(0,-\gamma)$ is a ring homomorphism
$R\to\QQ[\gamma]$. It maps $G^{(-\alpha,-\beta)}_\rho$ to $G^{(0,\gamma)}_\rho$, so
by~\eqref{eq:duality} and~\eqref{eq:one-parameter-translation} it maps
$g^{(\alpha,\beta)}_\lambda$ to $\sg^{(\gamma)}_\lambda$, and hence
$b^{(k)}_{\lambda\mu}(\alpha,\beta)$ to the coefficient of $\sg^{(\gamma)}_\mu$ in
$p_k\sg^{(\gamma)}_\lambda$. Apply it to Proposition~\ref{prop:finite-rule} and
Theorem~\ref{thm:connected-rule}. A history containing a stay of weight
$-\alpha$ gets weight $0$. For a given word there is at most one history
without such stays: at each step the move is forced, namely the addition
if $\nu_i<\nu_{i-1}$ and the stay (of weight $\beta=-\gamma$) otherwise; this
is the deterministic history. By Lemma~\ref{lem:histories} it has
$k-|\mu/\lambda|$ stays, so its weight is $(-\gamma)^{k-|\mu/\lambda|}$. The
diagonal coefficient is the specialization of
Proposition~\ref{prop:diagonal}.
\end{proof}

For $\lambda=\one$ we obtain the expansion of the power sums in the one-parameter dual basis.

\begin{corollary}[cf.\ \cite{Lenart}]\label{cor:one-param-pk}
For every $k\geq1$,
\[
 p_k=\sum_{\substack{m\geq1,\ j\geq0\\ m+j\leq k}}
 (-1)^j(-\gamma)^{q}\binom{q+j}{j}\,\sg^{(\gamma)}_{(m,1^j)},
 \qquad q=k-m-j .
\]
Equivalently, as formal power series in $t$,
\[
 \sum_{k\geq1}p_kt^k=\sum_{m\geq1,\ j\geq0}(-1)^j\,
 \frac{t^{m+j}}{(1+\gamma t)^{j+1}}\,\sg^{(\gamma)}_{(m,1^j)} .
\]
\end{corollary}

\begin{proof}
Apply the substitution $(\alpha,\beta)\mapsto(0,-\gamma)$ of the previous
proof to Theorem~\ref{thm:empty-partition}. Then
$c^{(k)}_{m,j}$ becomes $(-1)^j[t^q](1+\gamma t)^{-j-1}=(-1)^j(-\gamma)^q\binom{q+j}{j}$
with $q=k-m-j$. Multiplying by $t^k$ and summing over $k$ gives the
generating function.
\end{proof}

For the standard normalization $\gamma=-1$ this gives
\[
 p_2=\sg_{(2)}-\sg_{(1,1)}+\sg_{(1)},\qquad
 p_3=\sg_{(3)}-\sg_{(2,1)}+\sg_{(1,1,1)}+\sg_{(2)}-2\,\sg_{(1,1)}+\sg_{(1)},
\]
where $\sg_\lambda=\sg^{(-1)}_\lambda$ is the dual stable Grothendieck function of
\cite{LamPylyavskyy}.

\begin{remark}
Corollary~\ref{cor:one-param-pk} also follows from Lenart's
Schur expansion \cite[Thm.~2.2]{Lenart} by Hall duality:
\[
 [\sg^{(\gamma)}_\mu]p_k
 =\sum_{a=0}^{k-1}(-1)^a
   [s_{(k-a,1^a)}]\,\sG^{(\gamma)}_\mu.
\]
Only hook partitions $\mu=(m,1^j)$ can contribute. If
$q=k-m-j\geq0$, the only contributing outer hook is
$(m,1^{j+q})$, and Lenart's flagged fillings give
\[
 [s_{(m,1^{j+q})}]\,\sG^{(\gamma)}_{(m,1^j)}
 =\gamma^q\binom{j+q}{q}.
\]
Multiplication by the hook sign $(-1)^{j+q}$ gives
$(-1)^j(-\gamma)^q\binom{j+q}{j}$, as claimed.
\end{remark}

\subsection{Grassmannian \texorpdfstring{$K$}{K}-homology}\label{subsec:K-homology}
At $(\alpha,\beta)=(0,1)$ we have $g^{(0,1)}_\lambda=\sg^{(-1)}_\lambda$. We use
the stable $K$-homological interpretation of Lam and Pylyavskyy
\cite[\S9.5, Thm.~9.16]{LamPylyavskyy}. For a rectangle $(n-m)^m$ and
$\lambda\subseteq(n-m)^m$, put
$\widetilde\lambda=(n-m-\lambda_m,\dots,n-m-\lambda_1)$. Their correspondence
sends $g^{(0,1)}_\lambda$ to the class of the ideal sheaf of the boundary of
the Schubert variety $X_{\widetilde\lambda}\subseteq\mathrm{Gr}(m,n)$, and
multiplication to pushforward along the direct-sum maps between
Grassmannians. Products below are understood in this stable sense, not as
tensor products in the $K$-theory of a fixed Grassmannian. We write
$\eta_\lambda$ and $\pi_k$ for the classes corresponding to
$g^{(0,1)}_\lambda$ and $p_k$.

\begin{corollary}\label{cor:K-homology}
For every partition $\lambda$ and every $k\geq1$,
\[
 \pi_k\cdot\eta_\lambda=-d(\lambda)\,\eta_\lambda
 +\sum_\mu(-1)^{r-c}H_B(k-D-c+1)\,\eta_\mu ,
\]
the sum over the partitions $\mu\supsetneq\lambda$ such that $\mu/\lambda$ is a
corner-connected broken ribbon, with the statistics of
Definition~\ref{def:corner-connected}. In particular the coefficients are
integers and sign-coherent: $(-1)^{r-c}[\eta_\mu]\,\pi_k\eta_\lambda\geq0$
for $\mu\neq\lambda$. For $\mu\neq\lambda$ this coefficient is also the signed
number $\sum(-1)^{\leg(w)}$ of hook words $w$ of length $k$ with interval
support whose deterministic history from $\lambda$ ends at $\mu$.
\end{corollary}

\begin{proof}
The coefficient of $\eta_\mu$ is $b^{(k)}_{\lambda\mu}(0,1)$. The formula is
\eqref{eq:dual-stable-binomial} together with Proposition~\ref{prop:diagonal},
and the sign statement holds because $H_B(N)\geq0$. The last statement is
Corollary~\ref{cor:one-param-finite} with $\gamma=-1$.
\end{proof}

\subsection{Deformed power sums}\label{sec:further}

For an element $u$ of a commutative $\QQ$-algebra, the \emph{deformed power
sums} are
\begin{equation}\label{eq:deformed-power-sum}
 p^{\langle u\rangle}_k=\sum_{i\geq1}\Bigl(\frac{x_i}{1-ux_i}\Bigr)^k
 =\sum_{j\geq0}\binom{k+j-1}{j}u^j\,p_{k+j}\qquad(k\geq1);
\end{equation}
the second expression follows from
$(1-ux)^{-k}=\sum_{j\geq0}\binom{k+j-1}{j}u^jx^j$. If $u$ lies in the
coefficient ring, then $p^{\langle u\rangle}_k$ lies in the completion, and
its homogeneous components of degree $<k$ vanish. For $u=-\alpha$ these are
the multipliers in Kundu's rule (Section~\ref{sec:lowering}). Their action
on the $g$-basis is a generating function of the raising coefficients,
which Theorem~\ref{thm:compact-boundary} evaluates as a derivative of a
rational function.

Let $u$ be a further indeterminate. By~\eqref{eq:deformed-power-sum},
$p^{\langle u\rangle}_k\in\Lam_R[[u]]$, the ring of formal power series in $u$
with coefficients in $\Lam_R$. Define
\begin{equation}\label{eq:deformed-b}
 \widehat b^{(k)}_{\lambda\mu}(u)
 =\sum_{j\geq0}\binom{k+j-1}{j}u^j\,b^{(k+j)}_{\lambda\mu}\in R[[u]] .
\end{equation}
To specialize $u$ to elements such as $-\alpha$, we work over the ring of formal
power series $\QQ[[\alpha,\beta]]$. Let $\mathfrak m=(\alpha,\beta)$ be its
maximal ideal (series with zero constant term), and write $\operatorname{ord}$
for the $\mathfrak m$-adic order (the lowest total degree of a nonzero monomial
$\alpha^i\beta^j$, with $\operatorname{ord}(0)=+\infty$).

We use the $\mathfrak m$-adic topology on $\QQ[[\alpha,\beta]]$.
A series $\sum_j a_j$ converges if and only if
$\operatorname{ord}(a_j)\to+\infty$; equivalently, modulo
each power of $\mathfrak m$, only finitely many summands
are nonzero. Put
\[
 \hLam'=\prod_{D\geq0}\Lam^D_{\QQ[[\alpha,\beta]]}.
\]
Each homogeneous component is a finite free
$\QQ[[\alpha,\beta]]$-module in the Schur basis and carries
the coefficientwise $\mathfrak m$-adic topology. All
convergence assertions in $\hLam'$ below are understood
in this topology in each symmetric-function degree.

\begin{proposition}\label{prop:deformed}
Let $\lambda$ be a partition and $k\geq1$.
\begin{enumerate}[label=\textup{(\alph*)}]
\item In $\Lam_R[[u]]$,
\[
 p^{\langle u\rangle}_k\,g_\lambda=\sum_\mu\widehat b^{(k)}_{\lambda\mu}(u)\,g_\mu ,
\]
where the coefficient of each power of $u$ on the right is a finite sum.
Moreover $\widehat b^{(k)}_{\lambda\mu}(u)=0$ unless $\lambda\subseteq\mu$, and
\[
 \widehat b^{(k)}_{\lambda\lambda}(u)
 =d(\lambda)\biggl(\frac{(-\alpha)^k}{(1+\alpha u)^k}-\frac{\beta^k}{(1-\beta u)^k}\biggr).
\]
\item Let $u_0\in\mathfrak m$. For every partition $\mu$, the series
$\widehat b^{(k)}_{\lambda\mu}(u_0)$ converges in
$\QQ[[\alpha,\beta]]$. It is zero if
$\lambda\not\subseteq\mu$, and for $\lambda\subseteq\mu$ it satisfies
\[
 \operatorname{ord}
 \bigl(\widehat b^{(k)}_{\lambda\mu}(u_0)\bigr)
 \geq \bigl|\,|\mu/\lambda|-k\,\bigr|,
\]
where $\operatorname{ord}(0)=+\infty$.
Moreover,
\[
 p^{\langle u_0\rangle}_k\,g_\lambda
 =\sum_\mu
   \widehat b^{(k)}_{\lambda\mu}(u_0)\,g_\mu
 \qquad\text{in }\hLam',
\]
where the sum converges in each degree.
\end{enumerate}
\end{proposition}

\begin{proof}
(a) By~\eqref{eq:deformed-power-sum}, the coefficient of $u^j$ on the left
is $\binom{k+j-1}{j}p_{k+j}g_\lambda
=\binom{k+j-1}{j}\sum_\mu b^{(k+j)}_{\lambda\mu}g_\mu$, a finite sum by
Corollary~\ref{cor:support}; this is the coefficient of $u^j$ on the right.
The support statement follows from Corollary~\ref{cor:support}. By
Proposition~\ref{prop:diagonal} and
$\sum_{j\geq0}\binom{k+j-1}{j}y^j=(1-y)^{-k}$,
\[
 \widehat b^{(k)}_{\lambda\lambda}(u)
 =d(\lambda)\sum_{j\geq0}\binom{k+j-1}{j}u^j\bigl((-\alpha)^{k+j}-\beta^{k+j}\bigr)
 =d(\lambda)\bigl((-\alpha)^k(1+\alpha u)^{-k}-\beta^k(1-\beta u)^{-k}\bigr).
\]

(b) If $\lambda\not\subseteq\mu$, every coefficient in
\eqref{eq:deformed-b} vanishes. Suppose that
$\lambda\subseteq\mu$, and put
\[
 D=|\mu/\lambda|,\qquad
 j_0=\max(0,D-k),\qquad
 a_j=\binom{k+j-1}{j}u_0^j b^{(k+j)}_{\lambda\mu}.
\]
By Corollary~\ref{cor:support}, $a_j=0$ for $j<j_0$.
For $j\geq j_0$, homogeneity and
$\operatorname{ord}(u_0)\geq1$ give
\[
 \operatorname{ord}(a_j)\geq k+2j-D.
\]
This tends to infinity, so $\sum_j a_j$ converges.
Moreover,
\[
 k+2j-D\geq k+2j_0-D=|D-k|.
\]
Since $\mathfrak m^{|D-k|}$ is closed, the same lower bound
holds for the order of the sum.

It remains to justify the asserted expansion in $\hLam'$.
For each $j$, the ordinary raising rule gives the finite
identity
\[
 \binom{k+j-1}{j}u_0^j p_{k+j}g_\lambda
 =\sum_{\mu}
   \binom{k+j-1}{j}u_0^j
   b^{(k+j)}_{\lambda\mu}g_\mu.
\]
The sum of the left-hand sides over $j\geq0$ is
$p_k^{\langle u_0\rangle}g_\lambda$; in each symmetric-function
degree only finitely many $j$ contribute.

Fix a partition $\nu$, and write
$g_\mu=\sum_\nu c_{\mu\nu}s_\nu$. For every nonzero
summand on the right, Lemma~\ref{lem:grading} gives
\begin{align*}
 \operatorname{ord}\!\left(
   u_0^j b^{(k+j)}_{\lambda\mu}c_{\mu\nu}
 \right)
 &\geq j+(k+j-|\mu/\lambda|)+(|\mu|-|\nu|)\\
 &=2j+k+|\lambda|-|\nu|.
\end{align*}
This bound tends to infinity independently of $\mu$.
For each fixed $j$, only finitely many $\mu$ occur.
Consequently, modulo any fixed power $\mathfrak m^d$, only
finitely many pairs $(j,\mu)$ contribute to the coefficient
of $s_\nu$. We may therefore sum over $j$ and regroup by
$\mu$. The resulting coefficient of $g_\mu$ is precisely
$\widehat b^{(k)}_{\lambda\mu}(u_0)$, proving the identity
and its convergence in $\hLam'$.
\end{proof}

The compact generating function gives a rational expression
for each off-diagonal deformed coefficient.

\begin{corollary}
\label{cor:deformed-rational}
Let $u$ be a formal indeterminate and let $\lambda\subsetneq\mu$.
If $\mu/\lambda$ is not a corner-connected broken ribbon, then
$\widehat b^{(k)}_{\lambda\mu}(u)=0$. Otherwise,
\[
 \widehat b^{(k)}_{\lambda\mu}(u)
 =\frac{(-1)^{r-c}(\alpha+\beta)^{c-1}}{(k-1)!}
 \left(\frac{d}{du}\right)^{k-1}
 \left(
 \frac{u^{D+c-2}}{(1+\alpha u)^A(1-\beta u)^B}
 \right).
\]
The derivative treats $\alpha,\beta$ as constants. Specialization
$u=u_0\in(\alpha,\beta)\QQ[[\alpha,\beta]]$ is made after
differentiation, in the topology of Proposition~\ref{prop:deformed}.
\end{corollary}

\begin{proof}
Since $b^{(m)}_{\lambda\mu}=0$ for $1\leq m<D$, multiplying
\eqref{eq:compact-generating} (with $z=u$) by $u^{D-1}$ and setting $m=D+e$
gives
\[
 \sum_{m\geq1}b^{(m)}_{\lambda\mu}u^{m-1}
 =\frac{(-1)^{r-c}(\alpha+\beta)^{c-1}u^{D+c-2}}
 {(1+\alpha u)^A(1-\beta u)^B}.
\]
Differentiating $k-1$ times with respect to $u$ and dividing by $(k-1)!$
multiplies $b^{(m)}_{\lambda\mu}$ by $\binom{m-1}{k-1}u^{m-k}$ and kills the
terms with $m<k$, yielding~\eqref{eq:deformed-b}.
\end{proof}

\section{Multiplication in the canonical basis and a lowering rule}
\label{sec:lowering} The raising rule of Theorem~\ref{thm:A} is complementary to multiplication
by $p_k$ in the canonical basis, equivalently to the lowering problem for
the dual basis. In this section we deduce the latter from the $G$-side
rules of Nguyen--Hiep--Son--Thuy \cite{NHST} and Kundu \cite{Kundu}, by a
triangular change of multipliers and Hall duality, and compare the action
of Kundu's deformed power sums on the two bases.

\subsection{Connected skew diagrams and northwest rims}

The $G$-side rules are phrased in terms of connected skew diagrams and
their northwest rims. We record the facts about them that we need.

\begin{lemma}\label{lem:connected}
Let $\lambda\subsetneq\mu$ and $S=\{i:\mu_i>\lambda_i\}$. Then $\mu/\lambda$
is connected if and only if $S$ is an interval $[p,q]$ and
$\mu_{i+1}\geq\lambda_i+1$ for $p\leq i<q$. In this case
\[
 \rows(\mu/\lambda)=q-p+1,
 \qquad
 \cols(\mu/\lambda)=\mu_p-\lambda_q .
\]
\end{lemma}

\begin{proof}
Row $i\in S$ of $\mu/\lambda$ is the horizontal segment of columns
$[\lambda_i+1,\mu_i]$, which is connected. Two boxes in different rows
share an edge only if they lie in consecutive rows $i,i+1$ and in the same
column. Since $\lambda_{i+1}\leq\lambda_i$ and $\mu_{i+1}\leq\mu_i$, the
segments of rows $i$ and $i+1$ meet in the columns $[\lambda_i+1,\mu_{i+1}]$,
which is nonempty if and only if $\mu_{i+1}\geq\lambda_i+1$. Hence
$\mu/\lambda$ is connected if and only if the graph on $S$ joining $i$ and
$i+1$ whenever $\mu_{i+1}\geq\lambda_i+1$ is connected, which is the stated
condition. In that case the union of the segments is an interval of
columns, since consecutive segments overlap. Its right end is
$\max_{i\in S}\mu_i=\mu_p$ and its left end is
$\min_{i\in S}(\lambda_i+1)=\lambda_q+1$.
\end{proof}

The coefficients in \cite{NHST,Kundu} are described through the following ribbon inside a skew diagram.

\begin{definition}\label{def:NW}
For $\lambda\subseteq\mu$ define $\nu=\nu(\lambda,\mu)$ by
\[
 \nu_i=\min(\mu_i,\lambda_{i-1}+1)\qquad(i\geq1),
\]
with $\lambda_0=\infty$, so that $\nu_1=\mu_1$. The \emph{northwest rim} of
$\mu/\lambda$ is $\NW(\mu/\lambda)=\nu/\lambda$.
\end{definition}

\begin{lemma}\label{lem:NW}
Let $\lambda\subseteq\mu$ and $\nu=\nu(\lambda,\mu)$.
\begin{enumerate}[label=\textup{(\alph*)}]
\item $\nu$ is a partition with $\lambda\subseteq\nu\subseteq\mu$, the
diagram $\nu/\lambda$ contains no $2\times2$ block, and $\rho\subseteq\nu$ for
every partition $\rho$ with $\lambda\subseteq\rho\subseteq\mu$ such that
$\rho/\lambda$ contains no $2\times2$ block.
\item If $\mu/\lambda$ is connected and nonempty, then $\nu/\lambda$ is a
ribbon with $\rows(\nu/\lambda)=\rows(\mu/\lambda)$ and
$\cols(\nu/\lambda)=\cols(\mu/\lambda)$. It is the largest ribbon of the form
$\rho/\lambda$ with $\lambda\subseteq\rho\subseteq\mu$, and
\[
 |\NW(\mu/\lambda)|=\rows(\mu/\lambda)+\cols(\mu/\lambda)-1\leq|\mu/\lambda| .
\]
\end{enumerate}
\end{lemma}

\begin{proof}
(a) Since $\lambda_i\leq\mu_i$ and $\lambda_i<\lambda_{i-1}+1$, we have
$\lambda_i\leq\nu_i\leq\mu_i$. Moreover
$\nu_{i+1}=\min(\mu_{i+1},\lambda_i+1)\leq\min(\mu_i,\lambda_{i-1}+1)=\nu_i$
for $i\geq1$, so $\nu$ is a partition. By construction
$\nu_{i+1}\leq\lambda_i+1$ for all $i\geq1$, so $\nu/\lambda$ contains no
$2\times2$ block. If $\rho$ is as in the statement,
then $\rho_i\leq\mu_i$, and $\rho_i\leq\lambda_{i-1}+1$ for $i\geq2$; hence $\rho_i\leq\nu_i$.

(b) Because $\lambda_{i-1}+1>\lambda_i$, we have $\nu_i>\lambda_i$ if and only
if $\mu_i>\lambda_i$. So $\nu/\lambda$ and $\mu/\lambda$ occupy the same set of
rows, which is an interval $[p,q]$. For
$p\leq i<q$ Lemma~\ref{lem:connected} gives $\mu_{i+1}\geq\lambda_i+1$, hence
$\nu_{i+1}=\lambda_i+1$. By Lemma~\ref{lem:connected} again, $\nu/\lambda$ is
connected, and by (a) it is a ribbon. Next, $\nu_p=\mu_p$: this is clear if
$p=1$, and if $p\geq2$ then $p-1\notin[p,q]$ gives
$\mu_p\leq\mu_{p-1}=\lambda_{p-1}<\lambda_{p-1}+1$. Hence
$\cols(\nu/\lambda)=\nu_p-\lambda_q=\mu_p-\lambda_q=\cols(\mu/\lambda)$. Using
$\nu_{i+1}=\lambda_i+1$ for $p\leq i<q$,
\[
 |\nu/\lambda|=\nu_p+\sum_{i=p}^{q-1}(\lambda_i+1)-\sum_{i=p}^{q}\lambda_i
 =(\nu_p-\lambda_q)+(q-p)=\cols(\mu/\lambda)+\rows(\mu/\lambda)-1 .
\]
Finally, every ribbon $\rho/\lambda$ with $\lambda\subseteq\rho\subseteq\mu$
satisfies $\rho\subseteq\nu$ by (a), and $|\nu/\lambda|\leq|\mu/\lambda|$ since
$\nu\subseteq\mu$.
\end{proof}

We use the convention $\binom{n}{j}=0$ unless $0\leq j\leq n$. For
partitions $\lambda\subseteq\mu$, an integer $m\geq1$ and an element $t$ of a
commutative ring, put
\begin{equation}\label{eq:kappa-def}
 \kappa^{(m)}_{\lambda\mu}(t)=
 \begin{cases}
 t^{D-m}(-1)^{m-s}\dbinom{r-1}{m-s},
 &\begin{gathered}
   \mu/\lambda\text{ is connected and nonempty},\\
   s\leq m\leq r+s-1,
  \end{gathered}\\[6pt]
 0,&\text{otherwise},
 \end{cases}
\end{equation}
where $D=|\mu/\lambda|$, $r=\rows(\mu/\lambda)$, and
$s=\cols(\mu/\lambda)$ in the first case. By
Lemma~\ref{lem:NW}(b), $r+s-1\leq D$, so every displayed power of
$t$ has a nonnegative exponent. In particular the definition is valid
for every element $t$ of a commutative ring, including $t=0$.

\begin{remark}\label{rem:NW-automatic}
In \cite[Thm.~1.1]{NHST} and \cite[Corollary~2]{Kundu} the sum runs over the
partitions $\nu\supseteq\lambda$ such that $\nu/\lambda$ is connected,
$\cols(\nu/\lambda)\leq k$, and the maximal ribbon along the northwest border
of $\nu/\lambda$ has at least $k$ boxes; the coefficient contains the
factor $\binom{\rows(\nu/\lambda)-1}{k-\cols(\nu/\lambda)}$. The maximal
ribbon in question is $\NW(\nu/\lambda)$, by Lemma~\ref{lem:NW}. By the
same lemma, the conditions $\cols\leq k\leq|\NW(\nu/\lambda)|$ are
equivalent to $0\leq k-\cols(\nu/\lambda)\leq\rows(\nu/\lambda)-1$, that is,
to the nonvanishing of the binomial coefficient. Hence the sum in those
theorems may be taken over all $\nu$ with $\nu/\lambda$ connected, and the
coefficients are exactly $\kappa^{(k)}_{\lambda\nu}$ of~\eqref{eq:kappa-def}.
\end{remark}

\subsection{The \texorpdfstring{$G$}{G}-side rule}
Recall the deformed power sums $p_k^{\langle u\rangle}$ from
\eqref{eq:deformed-power-sum}. We use Kundu's two-parameter rule
in the following form.

\begin{theorem}[Kundu {\cite[Corollary~2]{Kundu}}; {\cite[Thm.~1.1]{NHST}} for $u=0$]
\label{thm:Kundu}
Let $u,v$ be indeterminates. For every partition $\lambda$ and $k\geq1$,
\[
 p^{\langle u\rangle}_k\,G^{(u,v)}_\lambda
 =\sum_{\nu\supsetneq\lambda}\kappa^{(k)}_{\lambda\nu}(-u-v)\,G^{(u,v)}_\nu
\]
in $\hLam_{\QQ[u,v]}$, where the sum converges degreewise.
\end{theorem}

\begin{proof}
Kundu proves this identity in $n$ variables for every $n\geq1$, with the sum
over the $\nu$ with $\ell(\nu)\le n$ described in Remark~\ref{rem:NW-automatic}
\cite[Corollary~2]{Kundu}; by that remark these are exactly the $\nu$ with
$\ell(\nu)\le n$ and $\kappa^{(k)}_{\lambda\nu}\ne0$. We deduce the identity in
$\hLam_{\QQ[u,v]}$.

By~\ref{G1}, $G^{(u,v)}_\nu$ has no component of degree $<|\nu|$. Hence the
homogeneous component of degree $D$ of the right-hand side involves only the
finitely many $\nu$ with $|\nu|\le D$.

For each $n\geq1$, let $\pi_n\colon\hLam_{\QQ[u,v]}\to\QQ[u,v][[x_1,\dots,x_n]]$
be the specialization map defined by setting $x_{n+1}=x_{n+2}=\cdots=0$.
Since $\pi_n$ maps every homogeneous symmetric polynomial of degree $D$ to a
homogeneous polynomial of degree $D$ (or to $0$), it preserves the
homogeneous components: $\pi_n(F)_D=\pi_n(F_D)$ for every $F\in\hLam_{\QQ[u,v]}$
and every $D\geq0$. By stability,
$\pi_n(G^{(u,v)}_\nu)=G^{(u,v)}_\nu(x_1,\dots,x_n)$, which vanishes if
$\ell(\nu)>n$, and $\pi_n(p^{\langle u\rangle}_k)=\sum_{i\le n}(x_i/(1-ux_i))^k$.
Thus $\pi_n$ maps the two sides of the theorem to the two sides of Kundu's
$n$-variable identity, and hence maps their difference $\Delta$ to $0$.

Fix $D\geq0$ and put $n=\max(1,D)$. Let
$\Delta_D\in\Lam^D_{\QQ[u,v]}$ be the degree-$D$ component of
$\Delta$. Then
\[
 \pi_n(\Delta_D)=(\pi_n\Delta)_D=0.
\]
Since every partition of $D$ has at most $n$ parts, the
specializations $m_\rho(x_1,\ldots,x_n)$, with $|\rho|=D$,
are linearly independent. Hence $\pi_n$ is injective on
$\Lam^D_{\QQ[u,v]}$, so $\Delta_D=0$ for every $D\geq0$.
Therefore $\Delta=0$.
\end{proof}
\begin{remark}\label{rem:deformed-Kundu}
The case of interest is $u_0=-\alpha$, for which
$p^{\langle-\alpha\rangle}_k=p_k\bigl(x/(1+\alpha x)\bigr)$ is the multiplier in
Kundu's rule for $G_\lambda=G^{(-\alpha,-\beta)}_\lambda$. On the canonical
basis it acts by the coefficients $\kappa^{(k)}_{\lambda\nu}(\alpha+\beta)$,
which depend only on the skew diagram $\nu/\lambda$ and vanish for
$\nu=\lambda$ (Theorem~\ref{thm:Kundu} with $(u,v)=(-\alpha,-\beta)$). On the
dual canonical basis, by Proposition~\ref{prop:deformed}, its diagonal
coefficient is
\[
 \widehat b^{(k)}_{\lambda\lambda}(-\alpha)
 =d(\lambda)\biggl(\frac{(-\alpha)^k}{(1-\alpha^2)^k}
 -\frac{\beta^k}{(1+\alpha\beta)^k}\biggr)\in\QQ[[\alpha,\beta]],
\]
which depends on $\lambda$ and not only on the (empty) skew diagram.
\end{remark}

\subsection{The lowering rule}

Kundu's rule concerns the deformed power sums. Expanding the ordinary power
sums in terms of them gives a rule for multiplication by $p_k$ in the
canonical basis, and by Hall duality a lowering rule for the dual basis.

Define, for partitions $\lambda\subseteq\mu$ and $k\geq1$,
\begin{equation}\label{eq:L-def}
 L^{(k)}_{\lambda\mu}(\alpha,\beta)
 =\sum_{j\geq0}\binom{k+j-1}{j}\,\alpha^j\,
  \kappa^{(k+j)}_{\lambda\mu}(\alpha+\beta),
\end{equation}
and $L^{(k)}_{\lambda\mu}=0$ if $\lambda\not\subseteq\mu$. By
Lemma~\ref{lem:NW}(b), only the terms with
$\cols(\mu/\lambda)\leq k+j\leq\rows(\mu/\lambda)+\cols(\mu/\lambda)-1\leq|\mu/\lambda|$
are nonzero. So the sum is finite, and $L^{(k)}_{\lambda\mu}=0$ unless
$\mu/\lambda$ is connected and $|\mu/\lambda|\geq k$. When it is nonzero,
$L^{(k)}_{\lambda\mu}$ is homogeneous of degree $|\mu/\lambda|-k$ and depends only on $|\mu/\lambda|$, $\rows(\mu/\lambda)$ and $\cols(\mu/\lambda)$.

\begin{proposition}\label{prop:G-side}
For every partition $\lambda$ and every $k\geq1$,
\[
 p_k\,G_\lambda=\sum_{\nu\supsetneq\lambda}L^{(k)}_{\lambda\nu}(\alpha,\beta)\,G_\nu
\]
in $\hLam_R$, where the sum converges degreewise.
\end{proposition}

\begin{proof}
Inverting $y=x/(1+\alpha x)$ as $x=y/(1-\alpha y)$ gives
$x^k=\sum_{j\geq0}\binom{k+j-1}{j}\alpha^jy^{k+j}$; summing over the
variables,
\begin{equation}\label{eq:power-sum-inversion}
 p_k=\sum_{j\geq0}\binom{k+j-1}{j}\alpha^j\,p^{\langle-\alpha\rangle}_{k+j}
\end{equation}
in $\hLam_R$, the sum converging degreewise because
$p^{\langle-\alpha\rangle}_m$ has no component of degree $<m$. Multiply by
$G_\lambda$ and apply Theorem~\ref{thm:Kundu} with $(u,v)=(-\alpha,-\beta)$
to each term:
\begin{equation}
    \label{eq:AA}
    p_kG_\lambda=\sum_{j\geq0}\binom{k+j-1}{j}\alpha^j
 \sum_{\nu\supsetneq\lambda}\kappa^{(k+j)}_{\lambda\nu}(\alpha+\beta)\,G_\nu.
\end{equation}
Here, for each $\nu$ such that $\nu/\lambda$ is connected, the coefficient $\kappa^{(k+j)}_{\lambda\nu}\neq0$ only for indices $k+j$ satisfying $\cols(\nu/\lambda) \leq k+j\leq \rows(\nu/\lambda) + \cols(\nu/\lambda) -1$, and by~\ref{G1}, $s_{\mu}$ appears in the Schur expansion of $G_{\nu}$, only for partitions $\nu$ with
$|\nu|\leq |\mu|$, so the coefficient of $s_{\mu}$ on the right of Eq.~\eqref{eq:AA} is a finite sum. Thus rearranging the sum is legitimate, in other words,
\[
p_kG_\lambda =\sum_{\nu\supsetneq\lambda} L^{(k)}_{\lambda\nu}(\alpha,\beta)\,G_\nu.
\]
\end{proof}

For $\alpha=0$ we have $L^{(k)}_{\lambda\nu}=\kappa^{(k)}_{\lambda\nu}(\beta)$, and
Proposition~\ref{prop:G-side} is the rule of \cite[Thm.~1.1]{NHST}, in
infinitely many variables. Pairing Proposition~\ref{prop:G-side} with the dual basis gives the lowering rule.

\begin{theorem}\label{thm:lowering}
For every partition $\mu$ and every $k\geq1$,
\begin{equation}\label{eq:lowering}
 p_k^{\perp}\,g_\mu
 =\sum_{\lambda\subsetneq\mu}L^{(k)}_{\lambda\mu}(\alpha,\beta)\,g_\lambda .
\end{equation}
\end{theorem}

\begin{proof}
For every partition $\lambda$, Lemma~\ref{lem:hopf}(b),
Proposition~\ref{prop:G-side} and~\eqref{eq:duality} give
\[
 [g_\lambda]\,p_k^\perp g_\mu
 =\Hall{G_\lambda}{p_k^\perp g_\mu}
 =\Hall{p_kG_\lambda}{g_\mu}
 =L^{(k)}_{\lambda\mu}(\alpha,\beta),
\]
where the last equality holds because $\Hall{G_\nu}{g_\mu}=\delta_{\nu\mu}$
and only the finitely many terms with $|\nu|\leq|\mu|$ enter the pairing. Since $L^{(k)}_{\lambda\mu}=0$ unless $\lambda\subsetneq\mu$, this gives the
finite expansion~\eqref{eq:lowering}.
\end{proof}

\begin{example}\label{ex:lowering}
Take $\mu=(2,1)$ and $k=1$. The partitions $\lambda\subsetneq\mu$ with
$\mu/\lambda$ connected are $(2)$, $(1,1)$ and $\one$; the diagram
$(2,1)/(1)$ consists of two boxes meeting at a corner. For $\lambda=(2)$
and $(1,1)$ the diagram is a single box, and $L^{(1)}_{\lambda\mu}=1$. For
$\lambda=\one$ we have $|\mu/\lambda|=3$ and
$\rows(\mu/\lambda)=\cols(\mu/\lambda)=2$, so $\kappa^{(1)}=0$,
$\kappa^{(2)}=\alpha+\beta$ and $\kappa^{(3)}=-1$. Hence
\[
 L^{(1)}_{\one,(2,1)}=\alpha(\alpha+\beta)-\alpha^2=\alpha\beta,
 \qquad
 p_1^\perp g_{(2,1)}=g_{(2)}+g_{(1,1)}+\alpha\beta .
\]
\end{example}

\begin{remark}\label{rem:lowering-classical}
At $\alpha=\beta=0$ we have $L^{(k)}_{\lambda\mu}(0,0)=\kappa^{(k)}_{\lambda\mu}(0)$,
which is nonzero only if $|\mu/\lambda|=k$. In that case
$k\leq\rows(\mu/\lambda)+\cols(\mu/\lambda)-1=|\NW(\mu/\lambda)| \leq |\mu/\lambda|$ forces $\NW(\mu/\lambda)=\mu/\lambda$,
i.e.\ $\mu/\lambda$ is a ribbon, and then $k-\cols(\mu/\lambda)=\rows(\mu/\lambda)-1$ and
$L^{(k)}_{\lambda\mu}(0,0)=(-1)^{\hgt(\mu/\lambda)}$. This is the adjoint of
the classical rule~\eqref{eq:classical-MN}.
\end{remark}

For one parameter, Theorem~\ref{thm:lowering} is the Hall dual of the rule of \cite{NHST}.

\begin{corollary}
\label{cor:lowering-one}
With $\sg^{(\gamma)}_\mu$ as in~\eqref{eq:one-parameter-translation},
\[
 p_k^\perp\,\sg^{(\gamma)}_\mu
 =\sum_{\lambda}(-\gamma)^{|\mu/\lambda|-k}(-1)^{k-\cols(\mu/\lambda)}
 \binom{\rows(\mu/\lambda)-1}{k-\cols(\mu/\lambda)}\,\sg^{(\gamma)}_\lambda.
\]
The sum is over $\lambda\subsetneq\mu$ such that $\mu/\lambda$ is
connected and
\[
 \cols(\mu/\lambda)\leq k\leq
 \rows(\mu/\lambda)+\cols(\mu/\lambda)-1.
\] 
\end{corollary}

\begin{proof}
Specialize $(\alpha,\beta)=(0,-\gamma)$ in~\eqref{eq:lowering}. This
specialization is a ring homomorphism $R\to\QQ[\gamma]$ which maps
$g^{(\alpha,\beta)}_\mu$ to $\sg^{(\gamma)}_\mu$, since~\eqref{eq:duality}
specializes to the duality defining $\sg^{(\gamma)}$. In~\eqref{eq:L-def}
only the term $j=0$ survives, and it equals $\kappa^{(k)}_{\lambda\mu}(-\gamma)$.
\end{proof}

\begin{remark}\label{rem:refined}
The theorem of Iwao, Motegi and Scrimshaw holds for refined parameters
$\balpha=(\alpha_1,\alpha_2,\dots)$ and $\bbeta=(\beta_1,\beta_2,\dots)$
\cite{IMS25}. It would be interesting to extend
Theorems~\ref{thm:A} and~\ref{thm:B} to that setting.
\end{remark}

\appendix

\section{Small expansions}\label{app:tables}

The following expansions were computed from Theorem~\ref{thm:compact-boundary}
and Proposition~\ref{prop:diagonal}; the case $\lambda=\one$ is
Example~\ref{ex:pk-small}. It can be checked that
Proposition~\ref{prop:finite-rule}, Theorem~\ref{thm:connected-rule} and
Theorem~\ref{thm:boundary} give the same coefficients for all $|\lambda|\leq5$
and $k\leq5$. As an independent check, we expanded $G^{(-\alpha,-\beta)}_\rho$
in Schur functions from the bialternant
formula~\eqref{eq:canonical-determinant}, computed $g_\lambda$ by
inverting~\eqref{eq:duality}, computed $p_kg_\lambda$ and $p_k^\perp g_\mu$ with
the classical rule~\eqref{eq:classical-MN} and its adjoint, and re-expanded
in the $g$-basis. This agrees with Theorem~\ref{thm:boundary} for all
$|\lambda|+k\leq7$, and with Theorem~\ref{thm:lowering} for all $|\mu|\leq7$. Theorem~\ref{thm:compact-boundary} was also compared with the interval form
for all $|\lambda|\leq6$ and $k\leq7$.

{\small
\begin{align*}
p_1g_{1} &= g_{2}+g_{11}-(\alpha+\beta)\,g_{1}\\[2pt]
p_2g_{1} &= g_{3}-g_{111}-2\alpha\,g_{2}+2\beta\,g_{11}+(\alpha^{2}-\beta^{2})\,g_{1}\\[2pt]
p_3g_{1} &= g_{4}-g_{22}+g_{1111}-3\alpha\,g_{3}+(\alpha+\beta)\,g_{21}-3\beta\,g_{111}\\
 &\quad+3\alpha^{2}\,g_{2}+3\beta^{2}\,g_{11}-(\alpha^{3}+\beta^{3})\,g_{1}\\[2pt]
p_1g_{2} &= g_{3}+g_{21}-(\alpha+\beta)\,g_{2}\\[2pt]
p_2g_{2} &= g_{4}+g_{22}-g_{211}-2\alpha\,g_{3}-(\alpha-\beta)\,g_{21}+(\alpha^{2}-\beta^{2})\,g_{2}\\[2pt]
p_3g_{2} &= g_{5}-g_{221}+g_{2111}-3\alpha\,g_{4}-(\alpha-2\beta)\,g_{22}+(\alpha-2\beta)\,g_{211}\\
 &\quad+3\alpha^{2}\,g_{3}+(\alpha^{2}-\alpha\beta+\beta^{2})\,g_{21}-(\alpha^{3}+\beta^{3})\,g_{2}\\[2pt]
p_1g_{11} &= g_{21}+g_{111}-(\alpha+\beta)\,g_{11}\\[2pt]
p_2g_{11} &= g_{31}-g_{22}-g_{1111}-(\alpha-\beta)\,g_{21}+2\beta\,g_{111}+(\alpha^{2}-\beta^{2})\,g_{11}\\[2pt]
p_3g_{11} &= g_{41}-g_{32}+g_{11111}-(2\alpha-\beta)\,g_{31}+(2\alpha-\beta)\,g_{22}-3\beta\,g_{1111}\\
 &\quad+(\alpha^{2}-\alpha\beta+\beta^{2})\,g_{21}+3\beta^{2}\,g_{111}-(\alpha^{3}+\beta^{3})\,g_{11}\\[2pt]
p_1g_{21} &= g_{31}+g_{22}+g_{211}-2(\alpha+\beta)\,g_{21}\\[2pt]
p_2g_{21} &= g_{41}-g_{2111}-2\alpha\,g_{31}-(\alpha-\beta)\,g_{22}+2\beta\,g_{211}+2(\alpha^{2}-\beta^{2})\,g_{21}\\[2pt]
p_3g_{21} &= g_{51}-g_{33}-g_{222}+g_{21111}-3\alpha\,g_{41}+(\alpha+\beta)\,g_{32}+(\alpha+\beta)\,g_{221}\\
 &\quad-3\beta\,g_{2111}+3\alpha^{2}\,g_{31}+(\alpha^{2}-\alpha\beta+\beta^{2})\,g_{22}+3\beta^{2}\,g_{211}
 -2(\alpha^{3}+\beta^{3})\,g_{21}
\end{align*}
}

In each line the terms of top degree reproduce the classical
rule~\eqref{eq:classical-MN} (Corollary~\ref{cor:top}), the diagonal term is
$d(\lambda)((-\alpha)^k-\beta^k)$, and the lines for $\lambda=(1,1)$ are obtained
from those for $\lambda=(2)$ by Corollary~\ref{cor:omega}. The coefficient
$\alpha+\beta$ of $g_{32}$ in $p_3g_{21}$ is Example~\ref{ex:disconnected}.

\end{document}